\documentclass[10pt,twoside]{article}
\usepackage[latin1]{inputenc}
\usepackage{amsmath}
\usepackage{graphicx}
\usepackage{yhmath}
\usepackage[table]{xcolor}
\usepackage{mathrsfs} 
\usepackage{amssymb}
\usepackage{url}
\usepackage{makecell}
\usepackage{arydshln}
\usepackage{multirow}
\usepackage{arydshln}
\usepackage{mathtools}
\usepackage{stmaryrd}  

\input epsf
\def\limiten{\renewcommand{\arraystretch}{0.5}
\begin{array}[t]{c}\stackrel{}{\longrightarrow} \\
{\scriptstyle n\rightarrow
\infty}\end{array}\renewcommand{\arraystretch}{1}}

\def\limitepsn{\renewcommand{\arraystretch}{0.5}
\begin{array}[t]{c}\stackrel{a.s.}{\longrightarrow} \\
{\scriptstyle n \rightarrow
\infty}\end{array}\renewcommand{\arraystretch}{1}}

\def\limiteloin{\renewcommand{\arraystretch}{0.5}
\begin{array}[t]{c}\stackrel{{\cal D}}{\longrightarrow} \\
{\scriptstyle n\rightarrow
\infty}\end{array}\renewcommand{\arraystretch}{1}}

\def\limiteproban{\renewcommand{\arraystretch}{0.5}
\begin{array}[t]{c}\stackrel{{\cal P}}{\longrightarrow} \\
{\scriptstyle n\rightarrow
\infty}\end{array}\renewcommand{\arraystretch}{1}}

\def\limitem{\renewcommand{\arraystretch}{0.5}
\begin{array}[t]{c}\stackrel{}{\longrightarrow} \\
{\scriptstyle m\rightarrow
\infty}\end{array}\renewcommand{\arraystretch}{1}}

\def\limitepsm{\renewcommand{\arraystretch}{0.5}
\begin{array}[t]{c}\stackrel{a.s.}{\longrightarrow} \\
{\scriptstyle m \rightarrow
\infty}\end{array}\renewcommand{\arraystretch}{1}}

\def\limiteloim{\renewcommand{\arraystretch}{0.5}
\begin{array}[t]{c}\stackrel{{\cal D}}{\longrightarrow} \\
{\scriptstyle m\rightarrow
\infty}\end{array}\renewcommand{\arraystretch}{1}}

\def\limiteloimds{\renewcommand{\arraystretch}{0.5}
\begin{array}[t]{c}\stackrel{{\cal D(S)}}{\longrightarrow} \\
{\scriptstyle m\rightarrow
\infty}\end{array}\renewcommand{\arraystretch}{1}}

\numberwithin{equation}{section}

\newtheorem{thm}{Theorem}[section]

\newtheorem{Corol}[thm]{Corollary}

\newtheorem{lem}[thm]{Lemma}

\newtheorem{rmrk}[thm]{Remark}

\newcommand{\E}{\ensuremath{\mathbb{E}}}
\newcommand{\R}{\ensuremath{\mathbb{R}}}
\newcommand{\Z}{\ensuremath{\mathbb{Z}}}

\newcommand{\N}{\ensuremath{\mathbb{N}}}

\newcommand{\var}{\ensuremath{\mathrm{Var}}}
\newcommand{\prob}{\ensuremath{\mathbb{P}}}

\definecolor{grisclair}{gray}{0.9}

\renewcommand{\arraystretch}{.8}

\begin{document}
\title{\bf Poisson QMLE for change-point detection in general integer-valued time series models}
 \maketitle \vspace{-1.0cm}
\begin{center}
   Mamadou Lamine DIOP \footnote{Supported by
   the MME-DII center of excellence (ANR-11-LABEX-0023-01) 
   } 
   and 
    William KENGNE \footnote{Developed within the ANR BREAKRISK: ANR-17-CE26-0001-01} 
 \end{center}

  \begin{center}
  { \it 
 THEMA, CY Cergy Paris Université, 33 Boulevard du Port, 95011 Cergy-Pontoise Cedex, France\\
  E-mail: mamadou-lamine.diop@u-cergy.fr ; william.kengne@u-cergy.fr  \\
  }
\end{center}

 \pagestyle{myheadings}
 \markboth{PQMLE for change-point detection in general integer-valued time series models}{Diop and Kengne}

~~\\
\textbf{Abstract}:
  We consider together the retrospective and the sequential change-point detection in a general class of integer-valued time series.
 The conditional mean of the process depends on a parameter $\theta^*$ which may change over time. We propose procedures which are based on the Poisson quasi-maximum likelihood estimator of the parameter, and where the updated estimator is computed without the historical observations in the sequential framework. For both the retrospective and the sequential detection, the test statistics converge to some distributions obtained from the standard Brownian motion under the null hypothesis of no change and diverge to infinity under the alternative; that is, these procedures are consistent.
 Some results of simulations as well as real data application are provided.
 
 \medskip
 

\section{Introduction}
%
We consider a class of integer-valued time series in a semiparametric framework.
Let $\Theta$ be a fixed compact subset of $\R^d$ ($d \in \N$) and $\mathcal T \subseteq \Z$. 
 For any $\theta \in \Theta$, define the class of observation-driven models given by
 
 \medskip
 {\bf Class} $\mathcal{OD}_{\mathcal T }(f_{\theta})$: A $\N_0$-valued ($\N_0=\N \cup \{0\}$) process $Y=\{Y_{t},\,t\in \Z \}$ belongs to $\mathcal{OD}_{\mathcal T }(f_{\theta})$ if it satisfies:
  \begin{equation} \label{Model}
             \E(Y_t|\mathcal{F}_{t-1})=f_{\theta}(Y_{t-1},Y_{t-2},\cdots) ~~ \forall t \in \mathcal{T} ,
   \end{equation}
where $\mathcal{F}_{t-1}=\sigma\left\{Y_{t-1},Y_{t-2},\cdots \right\}$ is the $\sigma$-field generated by the whole past at time $t-1$,   
and $f_{\theta}(\cdot)$ is a measurable non-negative function, assumed to be known up to the  parameter $\theta$.  
This class includes numerous classical integer-valued time series, which can be written as a model $\mathcal{OD}_{\Z }(f_{\theta})$: for instance, Poisson, negative binomial, binary INGARCH or Poisson exponential autoregressive model (proposed by Fokianos {\it et al.} (2009)).  
The class $\mathcal{OD}_{\Z }(f_{\theta})$ has been studied by Ahmad and Francq (2016). 
 Under certain regularity conditions, they have established the consistency and asymptotic normality of the Poisson quasi-maximum likelihood estimator (PQMLE) of the model's parameter.\
 
 \medskip
 
In this work, our main focus of interest is the structural change-point problem in the model (\ref{Model}). 
 By relying on the PQMLE, we will address this issue in two different frameworks: the retrospective (or off-line) and sequential (or on-line) detection. The retrospective detection is performed when all the data are available, whereas  the on-line approach focus on sequential change detection as long as new data arrive. For surveys on these approaches, we refer readers to Basseville and Nikiforov (1993) and Cs{\"o}rg{\"o} and Horv{\'a}th (1997).
 
  \medskip
  
  The change-point problem in time series of count has been addressed in several studies; see  among others, Franke \textit{et al.} (2012),  Hudecov{\'a} (2013), Fokianos {\it et al.} (2014), Kang and Lee (2014),  Kirch and Tadjuidje Kamgaing (2014),  Doukhan and Kengne (2015), Cleynen and Lebarbier (2017), Diop and Kengne (2017, 2021) for some papers in the retrospective setting; and Kengne (2015), Kirch and Tajduidje Kamgaing (2015), Kirch and Weber (2018), Kengne and Ngongo (2020) for some recent papers in the sequential framework.
  Most of these works are developed in the parametric setting by assuming that the conditional distribution of the observation given the whole past is known; which is quite restrictive in practice. Diop and Kengne (2021) have considered the semiparametric framework, but they focussed on the model selection approach. Kirch and Tajduidje Kamgaing (2015) and Kirch and Weber (2018) developed a general setup based on estimating functions for sequential change-point detection in continuous and integer valued time series. 
  As pointed out by Kengne and Ngongo (2020), the optimal estimating function in several classical parametric model is based on the score function and in the case of infinite memory process considered here, a more complex class of estimating functions is needed; which can lead  some difficulties in the application of such sequential procedure.
  
  \medskip
 In this contribution, we consider a process $Y=\{Y_{t},\,t\in \Z \}$ satisfying (\ref{Model}) depending on a parameter $\theta^*$ which may change over time.
 \begin{itemize}
  \item[(i)] In the retrospective detection, we construct a statistics based on the PQMLE and establish that it converges to a well-known distribution under the null hypothesis (no change) and diverges to infinity under the alternative.
  \item[(ii)] In the sequential detection, we construct a detector, based on the PQMLE, which converges (to some distribution) under the null hypothesis and diverges to infinity under the alternative. In order to perform a procedure with a more efficient detection delay (see Theorem \ref{th4}), the updated estimator is computed without the historical observations.  
  \end{itemize} 
 For the both retrospective and sequential detection, the proposed procedure is consistent.
  
 \medskip

 The paper is structured as follows. Section 2 contains some classical assumptions as well as the definition of the PQMLE. 
 In Section 3, we derive the procedure for the retrospective change-point detection and provide the main results. 
 Section 4 focuses on the sequential change-point detection. 
 Some results of simulations and real data example are displayed in Section 5 whereas Section 6 is devoted to a concluding remarks. Section 7 provides the proofs of the main results. 
   

\section{Assumptions and Poisson QMLE}\label{Sect_Ass_PQMLE}
Throughout the sequel, the following notations will be used:
{\em
\begin{itemize}
 \item $ \|x \| \coloneqq  \sqrt{\sum_{i=1}^{p} |x_i|^2 } $, for any $x \in \mathbb{R}^{p}$;

\item  $\left\|f\right\|_{\Theta} \coloneqq \sup_{\theta \in \Theta}\left(\left\|f(\theta)\right\|\right)$ for any function $f:\Theta \longrightarrow   M_{p,q}(\R)$,
 where $M_{p,q}(\R)$ denotes the set of matrices of dimension $p\times q$ with coefficients in $\R$;

\item $\left\|Y\right\|_r \coloneqq \E\left(\left\|Y\right\|^r\right)^{1/r}$, where $Y$ is a random vector with finite $r-$order moments; 
\item $T_{\ell,\ell'}=\{\ell,\ell+1,\cdots,\ell'\}$ for any $(\ell,\ell') \in \N^2$ such as $\ell \leq \ell'$. 
\end{itemize}
}
\noindent 
We set the following classical Lipschitz-type condition on the function $f_\theta$.

    \medskip
    \noindent 
    \textbf{Assumption} \textbf{A}$_i (\Theta)$ ($i=0,1,2$):
    For any $y \in \mathbb{N}_0^{\N}$, the function $\theta \mapsto f_\theta(y)$ is $i$ times continuously differentiable on $\Theta$  with $ \left\| \partial^i f_\theta(0)/ \partial \theta^i\right\|_\Theta<\infty $; 
    and
      there exists a sequence of non-negative real numbers $(\alpha^{(i)}_k)_{k\geq 1} $ satisfying
     $ \sum\limits_{k=1}^{\infty} \alpha^{(0)}_k <1 $ (or $ \sum\limits_{k=1}^{\infty} \alpha^{(i)}_k <\infty $ for $i=1, 2$);
   such that for any  $y, y' \in \mathbb{N}_0^{\N}$,
  \[ 
  \sup_{\theta \in \Theta  } \Big \| \frac{\partial^i f_\theta(y)}{ \partial \theta^i}-\frac{\partial^i f_\theta(y')}{\partial\theta^i} \Big \|
  \leq  \sum\limits_{k=1}^{\infty}\alpha^{(i)}_k |y_k-y'_k| ,
  \]
where $\| \cdot\|$ denotes any vector, matrix norm.

 \medskip 
  \noindent 
In the whole paper, it is assumed that any $\{Y_{t} ,\,t\in \Z\}$ belonging to $\mathcal{OD}_{\mathcal T }(f_{\theta})$ is a stationary and ergodic process satisfying: 
 \begin{equation}\label{moment}
    \exists C>0, \epsilon >1, \text{ such that } \forall t \in \Z, ~ ~ \E Y_{t}^{1+\epsilon} <C. 
   \end{equation}
 
 \medskip 
 
 \noindent 
Let $k \geq 1$ and $\theta^* \in \Theta$. 
If $(Y_{1},\ldots,Y_{k}) \in \mathcal{OD}_{\{1,\cdots,k \}}(f_{\theta^{*}})$, then for any subset $\mathcal T \subseteq \{1,\cdots,k\}$, the conditional Poisson (quasi)log-likelihood computed on $\mathcal T$ is given (up to a constant) by
 \[
  L(\mathcal T, \theta) \coloneqq 
   \sum_{t \in \mathcal T} \ell_t(\theta) ~ \text{ with }~
  \ell_t(\theta) = Y_t\log \lambda_t(\theta)- \lambda_t(\theta),
  \]
  where $ \lambda_t(\theta) =f_\theta(Y_{t-1}, Y_{t-2}, \cdots )$.
 We approximate this conditional (quasi)log-likelihood (see Ahmad and Francq (2016), for more details) by
\begin{equation}\label{logvm}
\widehat{L}(\mathcal T, \theta) \coloneqq  \sum_{t\in \mathcal T} \widehat{\ell}_t(\theta) ~ \text{ with }~
 \widehat{\ell}_t(\theta) =  Y_t\log \widehat{\lambda}_t(\theta)- \widehat{\lambda}_t(\theta),
 \end{equation}
  where $ \widehat{\lambda}_t(\theta) =  f_\theta(Y_{t-1}, \cdots Y_{1},0,\cdots,0)$. 
 According to (\ref{logvm}), the PQMLE of $ \theta^*$ computed on $\mathcal T$ is defined by
 \begin{equation}\label{emv}
  \widehat{\theta}(\mathcal T) \coloneqq  \underset{\theta\in \Theta}{\text{argmax}} \big(\widehat{L}(\mathcal T,\theta)\big).
  \end{equation}

  \medskip
  
  \noindent
  If $(Y_{1},\ldots,Y_{n})$ is an observed trajectory of a process $\{Y_t,\, t \in \Z\}$ belonging to $\mathcal{OD}_{\Z}(f_{\theta^{*}})$, 
  then we set the following regularity assumptions to obtain the asymptotic results (consistency and asymptotic normality) of the PQMLE (see Ahmad and Francq (2016)):
  %
  \begin{enumerate}
    
    \item [(\textbf{A0}):] 
     for all  $(\theta, \theta')\in \Theta^2$,
 $ \big( f_\theta(Y_{t-1}, Y_{t-2}, \cdots)= f_{\theta'}(Y_{t-1}, Y_{t-2}, \cdots)  \ \text{a.s.} ~ \text{ for some } t \in \N \big) \Rightarrow ~ \theta = \theta'$; 
 moreover, $\exists  \underline{c}>0$ such that $\displaystyle \inf_{ \theta \in \Theta} f_\theta(y)  \geq \underline{c}$, for all $ y \in  \N_0^{\N} $;

\item [(\textbf{A1}):] $\theta^* $ is an interior point of $\Theta \subset \mathbb{R}^{d}$;
 
  \item [(\textbf{A2}):]  $a_{t} \overset{a.s}{\longrightarrow}  0$ and $Y_{t} a_{t} \overset{a.s}{\longrightarrow} 0$   as $t\rightarrow \infty$, where $a_{t}=\underset{\theta \in \Theta }{\sup} \left| \widehat{\lambda}_{t}(\theta) -\lambda_{t}(\theta)\right|$;
     \item [(\textbf{A3}):] The matrices $J = \E \Big[ \frac{1}{\lambda_{t}(\theta^* )}  \frac{\partial \lambda_{t}(\theta^* )}{ \partial \theta} \frac{\partial \lambda_{t}(\theta^* )}{ \partial \theta'}  \Big] $
           ~and~ $I = \E \Big[ \frac{\var(Y_{t}|\mathcal{F}_{t-1})}{\lambda^2_{t}(\theta^* )}  \frac{\partial \lambda_{t}(\theta^* )}{ \partial \theta} \frac{\partial \lambda_{t}(\theta^* )}{ \partial \theta'}  \Big]$ exist;
     \item [(\textbf{A4}):] for all $c' \in \R$, $c' \frac{\partial \lambda_{t} (\theta^* )}{\partial    \theta}=0$ a.s   $\Rightarrow ~ c'=0$;
     \item [(\textbf{A5}):] there exists a neighborhood $V(\theta^* )$ of $\theta^* $ such that: for all $i, k \in \left\{1,\cdots,d\right\} $, \[\E \left[ \sup_{\theta \in V(\theta^* )} \left|  \frac{\partial^2 }{ \partial \theta_i \partial \theta_k } \ell_{t}(\theta) \right|\right]<\infty ; \]
  %
    \item [(\textbf{A6}):]   $b_{t}$, $b_{t} Y_{t}$ and $a_{t}d_{t} Y_{t}$ are of order $O(t^{-h})$ for some $h>1/2$, where
    \[
    b_{t}=\underset{\theta \in \Theta }{\sup} \left\{\E \left[\left\|\frac{\partial \widehat{\lambda}_{t} (\theta)}{ \partial \theta} -\frac{\partial \lambda_{t}(\theta)}{ \partial \theta}  \right\| \right]\right\}~~\textrm{and}~~d_{t}=\underset{\theta \in \Theta }{\sup}\max \left\{ \E\left[\left\| \frac{1}{\widehat{\lambda}_{t}(\theta)} \frac{\partial \widehat{\lambda}_{t} (\theta)}{ \partial \theta} \right\|\right],  \E\left[ \left\| \frac{1}{\lambda_{t}(\theta)} \frac{\partial \lambda_{t} (\theta)}{ \partial \theta} \right\|\right] \right\}. 
     \]
    
    \end{enumerate}
    
    \medskip
\begin{rmrk}
   The aforementioned assumptions have been imposed by Ahmad and Francq (2016) to study the asymptotic behavior of the PQMLE; in their works, they proved that all assumptions hold for many classical models. But, many of these assumptions (more precisely,  (\textbf{A2}), (\textbf{A3}), (\textbf{A5}) and (\textbf{A6})) can be easily obtained from the Lipschitz-type condition  \textbf{A}$_i (\Theta)$ with $i=0,1,2$. 
 \end{rmrk}
 Under (\textbf{A0})-(\textbf{A6}), Ahmad and Francq (2016) have established that the estimator  $\widehat{\theta}(T_{1,n})$ is strongly consistent and asymptotically normal; that is,
\begin{equation}\label{res_Francq}
 \widehat{\theta}(T_{1,n}) \limitepsn  \theta^* ~ 
\text{ and }~
 \sqrt{n}(\widehat{\theta}(T_{1,n})-\theta^*) \limiteloin \mathcal{N}(0,\Sigma^{-1}), \text{ with } \Sigma \coloneqq J I^{-1} J; 
 \end{equation}
where $I$ and $J$ are defined in the assumption (\textbf{A3}).   
 According to (\textbf{A4}), one can show that the matrices $I$ and $J$ are symmetric and positive definite.  
Throughout the sequel, we set for any $\ell, \ell' \in \N$ with $\ell \leq \ell'$, 
\begin{align*}
   &
   \widehat J(T_{\ell,\ell'})  =  \frac{1}{\ell'-\ell+1}\sum_{t \in T_{\ell,\ell'}}\frac{1}{\widehat \lambda_{t}(\widehat{\theta}(T_{\ell,\ell'}))}  \frac{\partial \widehat \lambda_{t}(\widehat{\theta}(T_{\ell,\ell'}))}{ \partial \theta} \frac{\partial \widehat \lambda_{t}(\widehat{\theta}(T_{\ell,\ell'}))}{ \partial \theta'} ,\\
   &
  \widehat I(T_{\ell,\ell'})  = \frac{1}{\ell'-\ell+1}\sum_{t \in T_{\ell,\ell'}} \Big(\frac{Y_t}{\widehat \lambda_{t}(\widehat{\theta}(T_{\ell,\ell'}))} -1 \Big)^2 \,  \frac{\partial \widehat \lambda_{t}(\widehat{\theta}(T_{\ell,\ell'}))}{ \partial \theta} \frac{\partial \widehat \lambda_{t}(\widehat{\theta}(T_{\ell,\ell'}))}{ \partial \theta'} .
 \end{align*}
Under the previous assumptions, $\widehat J(T_{1,n})$ and  $\widehat I(T_{1,n})$ converge almost surely to $J$ and $I$ respectively.
 Hence, the matrix $\Sigma $ can be consistently estimated by $\widehat \Sigma_n  = \widehat J(T_{1,n}) \widehat I(T_{1,n})^{-1}   \widehat J(T_{1,n})$.


\section{Poisson QMLE for retrospective change-point detection}
Assume that a trajectory $(Y_{1},\cdots,Y_{n})$ of the process $\{Y_{t},\,t\in \Z \}$ is observed and consider the following test problem:
\begin{enumerate}
    \item [ H$_0$:] $(Y_1,\cdots,Y_n)$  is a trajectory the  process $\{Y_{t},\,t\in \Z \} \in \mathcal{OD}_{\Z}(f_{\theta_1^*})$  with $\theta_1^* \in \Theta$.
    
    \item [ H$_1$:] There exists $((\theta^{*}_1,\theta^{*}_2),t^{*}) \in \Theta^{2}\times \{2,3,\cdots, n-1 \}$ (with $\theta^{*}_1 \neq \theta^{*}_2$) such that
     $(Y_1,\cdots,Y_{t^{*}})$ is a trajectory of a process $\{Y^{(1)}_{t},\, t \in \Z\} \in \mathcal{OD}_{\Z}(f_{\theta^*_1})$   and  $(Y_{t^{*}+1},\cdots,Y_n)$ a trajectory of a process 
      $\{Y^{(2)}_{t}, \,t \in \Z\} \in \mathcal{OD}_{\Z}(f_{\theta^*_2})$. 
\end{enumerate}
 By using the PQLME of the parameter, we construct a semi-parametric test statistic from the basic idea that, under the null hypothesis (i.e., no change), for $1<k<n$, $\widehat{\theta}(T_{1,k})$ and  $\widehat{\theta}(T_{k+1,n})$ are close to $\widehat{\theta}(T_{1,n})$; 
 that is, the distances $\|\widehat{\theta}(T_{1,k})-\widehat{\theta}(T_{1,n}) \|$ and $\|\widehat{\theta}(T_{k+1,n})-\widehat{\theta}(T_{1,n}) \|$ are expected to not be too large. 

\medskip
\noindent
Let $(u_n)_{n\geq 1}$ and $(v_n)_{n\geq 1}$ be two integer valued sequences satisfying: $u_n,v_n \rightarrow +\infty$ and $\frac{u_n}{n},\frac{v_n}{n} \rightarrow 0$ as $n\rightarrow +\infty$.
For all $n \geq 1$, define the matrix $\widehat{\Sigma}(u_n)$  by
\begin{equation}\label{Sigma_un}
\widehat{\Sigma}(u_n)=\frac{1}{2} 
\left[
\widehat J(T_{1,u_n})  \widehat I(T_{1,u_n})^{-1}   \widehat J(T_{1,u_n}) +
\widehat J(T_{u_n+1,n})  \widehat I(T_{u_n+1,n})^{-1}  \widehat J(T_{u_n+1,n})
\right].
\end{equation}

\medskip

\noindent
Now, consider the test statistic:
 \[ \widehat{C}_n=\max_{v_n\leq k \leq n-v_n}\widehat{C}_{n,k} \]
with
 \[ \widehat{C}_{n,k}=\frac{1}{q^{2}(\frac{k}{n})}\frac{k^{2}(n-k)^{2}}{n^{3}}\left(\widehat{\theta}(T_{1,k})-\widehat{\theta}(T_{k+1,n})\right)' \widehat{\Sigma}(u_n) \left(\widehat{\theta}(T_{1,k})-\widehat{\theta}(T_{k+1,n})\right),
  \]
 where $q:(0, 1) \mapsto (0, \infty)$ is a weight function non-decreasing in a neighborhood of zero, non-increasing in a neighborhood of one and satisfying
  \[ 
  \inf_{\varphi<\tau<1-\varphi}q(\tau)>0 ~~\textrm{for all}~~ 0<\varphi< \frac{1}{2}.
  \]
Its behavior can be controlled at the neighborhood of zero and one by the integral (see Cs{\"o}rg{\"o} \textit{et al.} (1986))
 \[
   I(q,c)=\int_{0}^{1}\frac{1}{t(1-t)}\exp\left(-\frac{cq^{2}(t)}{t(1-t)}\right)dt ,~~c>0.
 \]
The weight function $q$ allows to increase the power of the test procedure based on the statistic $\widehat{C}_n$.

 The proprieties of the matrix $\widehat{\Sigma}(u_n)$ is very important to prove of the consistency of the proposed procedure. Indeed, when the parameter is constant over the observations (under H$_0$), according to the assumptions of Section \ref{Sect_Ass_PQMLE}, we can show that $\widehat{\Sigma}(u_n)$ is also a consistent estimator of the covariance matrix $\Sigma$.  
Under the alternative, the model depends on two parameters and the consistency of $\widehat{\Sigma}(u_n)$ is not ensured. 
But, under the classical Assumption $\bf B$ (see below), one can show that the first matrix on the right hand side of (\ref{Sigma_un}) 
converges to the covariance matrix of the stationary model of the first regime which is positive definite and the second matrix is positive semi-definite. This will play a key role in the proof of the consistency under the alternative. \\
Let us note that the sequence $v_n$ is also very important for the statistic $\widehat{C}_n$;  
it is used to assure that the length of $T_ {1, v_n} $ and $T_ {v_n+1, n} $ are not too small, which allows to obtain the convergence of the  numerical algorithm used to compute these estimators.
Such approach to construct the statistic for change-point detection in a retrospective setting has already been used  by Doukhan and Kengne (2015) and Diop and Kengne (2017). 

\medskip
\noindent
Theorems \ref{th1} and \ref{th2} give the asymptotic behavior of the statistic $\widehat{C}_n$ under the null and alternative hypothesis.
\begin{thm}\label{th1}
Under H$_0$ with $\theta^*_1$  an interior point of $\Theta$, assume that (\textbf{A0})-(\textbf{A6}), \textbf{A}$_i(\Theta)$ ($i=0,1,2$) and (\ref{moment}) (with $\epsilon>2$) hold with
 \begin{equation}\label{eq_th1}
  \alpha_k^{(0)} + \alpha_k^{(1)} = \mathcal O (k^{-\gamma}), ~~ for~ some~ \gamma>3/2.
 \end{equation}
  If there exists $c>0$ such that $I(q,c)<\infty$, then
\[ \widehat{C}_n \limiteloin \sup_{0\leq\tau\leq1}\frac{\left\|W_d(\tau)\right\|^{2}}{q^{2}(\tau)},  \]
where $W_d$ is a $d$-dimensional Brownian bridge.
\end{thm}

\noindent 
According to the results of Theorem \ref{th1}, at a nominal level $\alpha \in (0,1)$, the critical region of the test is $(\widehat{C}_{n}>c_\alpha)$, where $c_\alpha$ is the $(1-\alpha)$-quantile of the distribution
of $\sup_{0<\tau<1} \big( \left\|W_d(\tau)\right\|^{2} / q^{2}(\tau) \big)$.  This assures that the test procedure has correct size asymptotically. 
In the empirical studies, we will consider the cases where $q\equiv 1$, and we will use the values of $c_\alpha$ provided in Lee \textit{et al.} (2003). 

\medskip
\noindent
Under the alternative, we assume

  \medskip
\noindent
 {\bf Assumption B}: {\em there exist $\tau \in (0,1)$ such that $t^*=[n\tau^*]$ (where $[x]$ is the integer part of $x$).}

\begin{thm}\label{th2}
Under H$_1$ with $\theta^*_1, \theta^*_2$ belong to the interior of $\Theta$, assume that {\bf B}, (\textbf{A0})-(\textbf{A6}),  \textbf{A}$_i(\Theta)$ ($i=0,1,2$), (\ref{moment}) (with $\epsilon>2$) and (\ref{eq_th1}) hold. Then,
\[ \widehat{C}_n  \limiteproban +\infty .\]
\end{thm}
This theorem establishes that the proposed procedure is consistent in power. Under H$_1$, a classical estimator of the breakpoint $t^*$ is
 \[   \widehat{t}_n =  \underset{v_n \leq k  \leq n - v_n }{ \text{argmax}  }  \widehat{C}_{n,k} . \]


\section{Sequential change-point detection}
Assume that we observed an available historical trajectory $(Y_{1},\ldots,Y_{m})$ generated from (\ref{Model}) according to a parameter $\theta^*_1$. New data $Y_{m+1},Y_{m+2},\cdots,$ will arrive and we would like to monitor these data from the time $m+1$ in order to test whether any structural change occurs. More precisely, for each new observation, we want to know if a change occurs in the parameter $\theta^*_1$. To address this problem, consider the following classical hypothesis testing:
\begin{enumerate}
    \item [ H$^*_0$:]  $\theta^*_1$ is constant over the observations $Y_1,\cdots,Y_{m},Y_{m+1},\cdots$; that is, $\{Y_{t},\,t\in \N \} \in \mathcal{OD}_{\N}(f_{\theta^*_1})$.
    \item [ H$^*_1$:] There exists $k^*>m$, $\theta^{*}_2 \in \Theta$ (with $\theta^{*}_1 \neq \theta^{*}_2$), such that
     $(Y_1,\cdots,Y_{k^*}) \in \mathcal{OD}_{\{1,\cdots,k^*\}}(f_{\theta^*_1})$ and  
     $\{Y_{k^*+n},\,n\in \N \} \in \mathcal{OD}_{\{k^*+1,\cdots\}}(f_{\theta^*_2})$. 
\end{enumerate}
Let $k>m$ be a monitoring instant.
As in Bardet and Kengne (2014), we derive a test procedure based on a statistic (called the detector) which evaluates the difference between $\widehat \theta(T_{\ell,k})$ and $\widehat \theta(T_{1,m})$ for any $\ell=m+1,\cdots,k$. 
Since the matrix $I(T_{1,m})$ and $J(T_{1,m})$ are symmetric and  non-singular (see Ahmad and Francq (2016)), from the central limit given in (\ref{res_Francq}), we deduce 
\[
 \sqrt{m}  \widehat I(T_{1,m})^{-1/2} \widehat J(T_{1,m})(\widehat{\theta}(T_{1,m})-\theta^*_0) \limiteloim \mathcal{N}(0,I_d), 
\] 
where $I_d$ is the identity matrix. Hence, define the detector
\[
 \widehat{D}_{k,\ell}= \sqrt{m}\frac{k-\ell}{k} \big \| \widehat I(T_{1,m})^{-1/2} \widehat J(T_{1,m})(\widehat{\theta}(T_{\ell,k})-\widehat{\theta}(T_{1,m})) \big \|, ~\text{ for all } k > m \text{ and } \ell =m+1,\cdots,k.
\]
To assure the convergence of $\widehat \theta(T_{\ell,k})$ and avoid some excessive distortion in the procedure, 
we introduce a sequence of integer numbers $(v'_m)_{m \geq 1}$ with $v'_m \ll m$ and define the set $\Pi_{m,k}=\left\{m-v'_m, m-v'_m+1,\cdots,k-v'_m\right\}$. Therefore, the detector $\widehat{D}_{k,\ell}$ will be computed for $\ell \in \Pi_{m,k}$. 
In the sequel, we assume that the sequence $(v'_m)_{m\geq 1}$ satisfies: 
\[
v'_m \rightarrow +\infty ~~~ \text{and}~~~  v'_m/\sqrt{m}
 \rightarrow 0 ~~~\text{as} ~~m\rightarrow +\infty.
\]
%
%
Let $T > 1$ ($T$ can be equal to infinity). The sequential monitoring scheme rejects  H$^*_0$ at the first time $k$ satisfying $m < k \leq [Tm]+1$ such that $\exists \ell \in\Pi_{m,k},\, \widehat{D}_{k,\ell}> c$ for a suitably chosen constant $c > 0$, where $[x]$ denote the integer part of $x$. The procedure is called \emph{closed-end} method when $T < \infty$ and
\emph{open-end} method when $T = \infty$. The set $\{m+1,m+2,\cdots, [Tm]\}$ is called the monitoring period, and its length depends on the time that we when to monitor the data.  
To obtain a procedure that is able to detect change that occurs at the beginning of the monitoring and that occurs a long time after the beginning of the monitoring, we will use a function $b : (0,\infty) \mapsto (0,\infty)$, called a boundary function satisfying:

\medskip
\noindent 
{\bf Assumption B$_*$}: {\em $b : (0,\infty) \mapsto (0,\infty)$ is a non-increasing and continuous function such as $\inf_{0<t<\infty} b(t)>0$.}

\medskip
\noindent
The monitoring scheme rejects H$^{*}_0$ at the first time $k$ (with $m < k \leq  [Tm] + 1$) such as there exists $\ell \in \Pi_{m,k}$ satisfying $\widehat{D}_{k,\ell}>  b((k-\ell)/m)$. Hence, we define the stopping time as follows: 
\begin{align*}
\tau(m) &=\inf \big\{ m< k \leq [Tm]+1 \, / \, \exists \ell \in \Pi_{m,k},\, \widehat{D}_{k,\ell}>  b((k-\ell)/m) \big \}\\
&=
\inf \big \{ m< k \leq [Tm]+1 \, / \, \max_{\ell \in \Pi_{m,k}} \frac{\widehat{D}_{k,\ell}}{b((k-\ell)/m)} >1  \big \},
\end{align*}
with the convention that $\inf \{ \emptyset \}=\infty$. 
 Therefore, we have
\begin{align} \label{Rejet_H0*}
\prob \{ \tau(m)< \infty \}
&= 
\prob \big\{ \max_{\ell \in \Pi_{m,k}} \frac{\widehat{D}_{k,\ell}}{b((k-\ell)/m)} >1  
\text{ for some } k \text{ between } m \text{ and }  [Tm]+1\big \} \nonumber\\
&=
\prob \bigg\{ \sup_{m<k\leq [Tm]+1}\max_{\ell \in \Pi_{m,k}} \frac{\widehat{D}_{k,\ell}}{b((k-\ell)/m)} >1 \bigg \}.
\end{align}
So, one would like to correctly calibrate a suitable boundary function $b(\cdot)$ such that the probability of false alarm is close to a fixed level $\alpha$ and the probability of true alarm is close to $1$, at least for $m$ large enough; that is, for some given $\alpha \in (0,1)$,
\begin{align} \label{cond_cons_Dkl}
\prob_{H^{*}_0} \{ \tau(m)< \infty \} \limitem \alpha 
~~\text{and}~~
\prob_{ H^{*}_1} \{ \tau(m)< \infty \} \limitem 1.
\end{align}
In the case where $b \equiv c$ with $c$ a positive constant, the first condition of (\ref{cond_cons_Dkl}) leads to compute the critical value $c=C_\alpha$ depending on $\alpha$. 
Moreover, if a change-point is detected under H$^{*}_1$; i.e., $\tau(m) < \infty$ and $\tau(m) > k^*$, then the detection
delay is defined by 
\begin{equation}\label{delay} 
\widehat d_m= \tau(m)-k^*.
\end{equation}
For the open and closed-end procedure, the following theorem gives the main result obtained under H$^{*}_0$
%
\begin{thm}\label{th3}
Under H$^*_0$  with $\theta^*_1$ an interior point of $\Theta$, assume that $\bf B_{*}$, (\textbf{A0})-(\textbf{A6}), \textbf{A}$_i(\Theta)$ ($i=0,1,2$),  (\ref{moment}) (with $\epsilon>2$) and (\ref{eq_th1}) hold.
 If $T = \infty$ (open-end procedure) or $T<\infty$ (closed-end procedure), then
	\begin{equation*}
	\prob \{ \tau(m)< \infty \} 
	\limitem 
	\prob \bigg\{ \sup_{1<t \leq T}\sup_{1<s<t} \frac{\| W_d(s)-sW_d(1)\|}{t b(s)} >1 \bigg \},
	\end{equation*}
	where $W_d$ is a d-dimensional standard Brownian motion.
\end{thm}
In the empirical studies, we will use the boundary function $b \equiv c$ with $c$ a positive constant. 
Since this function satisfies Assumption $\bf B_{*}$, the following corollary can be immediately deduced from Theorem \ref{th3}. 
 \begin{Corol}\label{corol1}
 Assume that $b(t)=c>0$, for all $t \geq0$. Under the assumptions of Theorem \ref{th3}, and with $T \in (1,\infty)$ or $T=\infty$,
 %
 \begin{equation*}
	\prob \{ \tau(m)< \infty \} 
	\limitem \prob \{ U_{d,T}> c \},
	\end{equation*}
	where, using the notation $U_{d,T}=U_{d,\infty}$ if  $T=\infty$,
	\begin{equation*}
	U_{d,T}= \sup_{1<t \leq T}\sup_{1<s<t} \frac{1}{t}\| W_d(s)-sW_d(1)\|.
	\end{equation*}
 \end{Corol}
 At a nominal level $\alpha \in (0, 1)$, we take $c = C_\alpha$, where $C_\alpha$ is the $(1-\alpha)$-quantile of the distribution of $U_{d,T}$. The values of $C_\alpha$ can be computed through Monte-Carlo simulations as described in Bardet and Kengne (2014). 
 
 \medskip
\noindent
 For the open-end and closed-end procedure, the following theorem shows that the detector tends to infinity when the parameter changes from $\theta^*_1$ to $\theta^*_2$ (under H$^*_1$).
 \begin{thm}\label{th4}
Under H$^*_1$ with $\theta^*_1, \theta^*_2$ belong to the interior of $\Theta$, assume that $\bf B_{*}$, (\textbf{A0})-(\textbf{A6}), \textbf{A}$_i(\Theta)$ ($i=0,1,2$), (\ref{moment}) (with $\epsilon>2$) and (\ref{eq_th1}) hold.
 If there exists $T^* \in(1,T)$ such that $k^*=k^*(m)=[T^*m]$, then for $k_m=k^*(m)+m^{\delta}$ wit $\delta \in (1/2,1)$,
	\begin{equation*}
	\max_{\ell \in \Pi_{m,k_m}} \frac{\widehat{D}_{k_m,\ell}}{b((k_m-\ell)/m)} 
	\limitepsm  \infty.
	\end{equation*}
\end{thm}
 The following corollary can be immediately deduced from the relation (\ref{Rejet_H0*}).
 \begin{Corol}\label{corol2}
 Under the assumptions of Theorem \ref{th4},
\[
	\prob \{ \tau(m)< \infty \} 
	\limitem 1.
	\]
	\end{Corol}
 Hence,  it follows from Theorem \ref{th4} that the change is detected with probability tending to one, both for open-end and closed-end (when $T^*<T$)  procedures and the detection delay $\widehat{d}_n$ can be bounded by
    $\mathcal{O}_P(m^{1/2+\varepsilon})$ for any $\varepsilon>0$ (or even by $\mathcal{O}_P\big (\sqrt m (\log m)^a \big )$ with $a>0$ using the same arguments).
	

 \section{Some numerical results}
In this section, we evaluate the performance of the proposed test procedures through an empirical study.
 For each procedure (the retrospective and the sequential detection),  we present some simulation results for change-point detection in a model that belongs to the class (\ref{Model}). 
 Applications to the number of transactions per minute for the stock Ericsson B  are also provided. The nominal level considered in the sequel is $\alpha = 0.05$.  
 %


\subsection{Simulation for the retrospective change-point detection}
The results of this subsection have been obtained by computing the test statistic $\widehat{C}_n$ with $q\equiv 1$ and $u_n, v_n$ equals to $[\left(\log(n)\right)^{\delta}]$ (with $2\leq \delta \leq 5/2$). 

\begin{figure}[h!]
\begin{center}
\includegraphics[height=6cm, width=17cm]{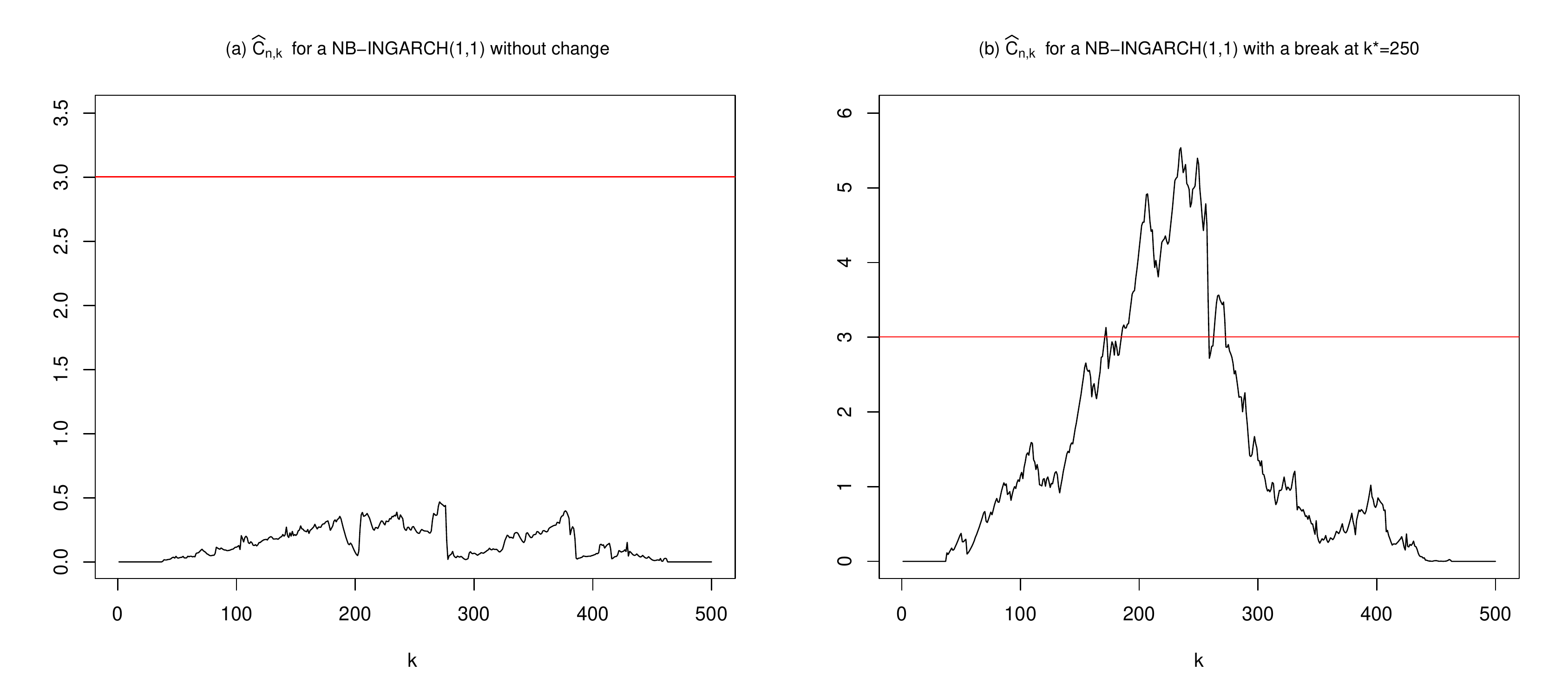}
\end{center}
\vspace{-.8cm}
\caption{\it \small Typical realization of the statistic $\widehat{C}_{n,k}$ for the retrospective change-point detection in  two NB-INGARCH$(1,1)$ processes with $r=1$.
\textbf{(a)} the corresponding values for a NB-INGARCH$(1,1)$ process without change, where the parameter
$\theta_0=(0.4,0.15,0.2)$ is constant. \textbf{(b)} the corresponding values for  a NB-INGARCH$(1,1)$ process where the parameter $\theta_0=(0.4,0.15,0.2)$ changes
to $\theta_1=(0.4,0.15,0.5)$ at $k^{*}=250$.
The horizontal line represents the limit of the critical region of the test.}
\label{Graphe_Offline_NB}
\end{figure}

\medskip
Consider the negative binomial INGARCH model (NB-INGARCH) given by 
 \begin{equation}\label{NB_INGARCH}
    Y_{t}|\mathcal{F}_{t-1}\sim\textrm{NB}(r,p_{t})~~\text{ with }  ~~ r\frac{(1-p_{t})}{p_{t}}=\lambda_{t}=\alpha^*_{0} +\alpha^* Y_{t-1} + \beta^* \lambda_{t-1},
    \end{equation}
    where $\alpha^*_{0}>0$, $\alpha^*, \beta^* \geq 0$ and 
 $NB(r,p)$ denotes the negative binomial distribution with parameters $r$ and $p$. We denote by $\theta = (\alpha^*_{0},\alpha^*,\beta^*)$ the parameter of the model. 
  We first generate two trajectories $(Y_1, \cdots, Y_{500})$ from (\ref{NB_INGARCH}): 
  a scenario under H$_0$ when the parameter is constant and a scenario under H$_1$ when the parameter changes at $k^{*} = 250$.
 The statistic $\widehat{C}_{n,k}$ is displayed in Figure \ref{Graphe_Offline_NB}. 
From this figure, one can see that,  
in the scenario without change, the statistic $\widehat{C}_{n}$ is less than the limit of the critical region that is represented by the horizontal line (see Figure \ref{Graphe_Offline_NB}(a)). 
Under the alternative (of change in the model), $\widehat{C}_{n}$ is greater than the critical value of the test and the statistic $\widehat{C}_{n,k}$ is large around the point where the change occurs (see Figure \ref{Graphe_Offline_NB}(b)). 

 \medskip
For $r=1,14$ and $n=500, 1000$, Table \ref{Table.1} indicates the empirical levels computed when the parameter is $\theta_0$ (under H$_0$) and the empirical powers computed when $\theta_0$ changes to $\theta_1$ at $t^{*} = 0.5n$ (under H$_1$); these results are based on $200$ replications. The scenario $\theta_0=(8.2, 0.2, 0.13)$; $\theta_1=(2.35, 0.12, 0.61)$ considered in the simulations  is close to the fitted model obtained from the number of transactions per minute for the stock Ericsson B (see below). 
 From the findings of Table \ref{Table.1}, one can see that even if the procedure exhibits some size distortion when $n=500$, the empirical levels are close to the nominal one when $n=1000$.  
 Also, the empirical powers displayed increases with the sample size.
 These results are consistent with Theorem \ref{th1}, \ref{th2} and are 
 overall satisfactory.
 %
 
\begin{table}[h!]
\centering
\scriptsize
\footnotesize
\caption{\small \it Empirical levels and powers at the nominal level $0.05$ for the retrospective change-point detection in a NB-INGARCH(1,1) process.}
\label{Table.1}
\vspace{.2cm}
\hspace*{-.3cm}
\begin{tabular}{lllcc}
\hline
\rule[0cm]{0cm}{.5cm}
&&$r$&$n=500$&$n=1000$\\
\hline
\rule[0cm]{0cm}{.4cm}
 Empirical levels: &&&&\\

                  &$\theta_0=(0.4,0.15,0.2) $ &$1$  &$0.065$&$0.045$\\
                  &                           &$14$ &$0.060$&$ 0.055$\\
                 
& & &&\\
                  & $\theta_0=(8.2, 0.2, 0.13)$ &$1$  &$0.080$&$0.070$\\
                  &                             &$14$ &$0.080$&$0.060$\\
                
&&&&\\
\rule[0cm]{0cm}{.1cm}
Empirical powers: & & &&\\
                   &$\theta_0=(0.4,0.15,0.2);~~~\theta_1=(0.4,0.15,0.5);$&$1$  &$0.605$&$0.885$\\
                   &                                                     &$14$ &$0.715$&$0.910$\\

& &&&\\

                  & $\theta_0=(8.2, 0.2, 0.13); ~~~\theta_1=(2.35, 0.12, 0.61);$&$1$  &$0.560$&$0.800$\\
                  &                                                             &$14$ &$0.870$&$0.955$\\

 \Xhline{.75pt}
\end{tabular}
\end{table}


\subsection{Simulation for the sequential change-point detection}
We consider the NB-INGARCH model and 
 focus on the  closed-end procedure with $T=1.5$; i.e, the historical available data are $Y_1,\cdots, Y_m$ and the monitoring period is $\left\{m+1,\cdots, [1.5m]\right\}$. The detector is computed with $v'_m = [(\log(m))^{\delta}]$ for $2\leq \delta \leq 5/2$.  
In the sequel, we denote  $\widehat{D}_{k} = \underset{\ell \in \Pi_{n,k}} {\max} \widehat{D}_{k,\ell} $, for any $k>m$.

\begin{figure}[h!]
\begin{center}
\includegraphics[height=6cm, width=17cm]{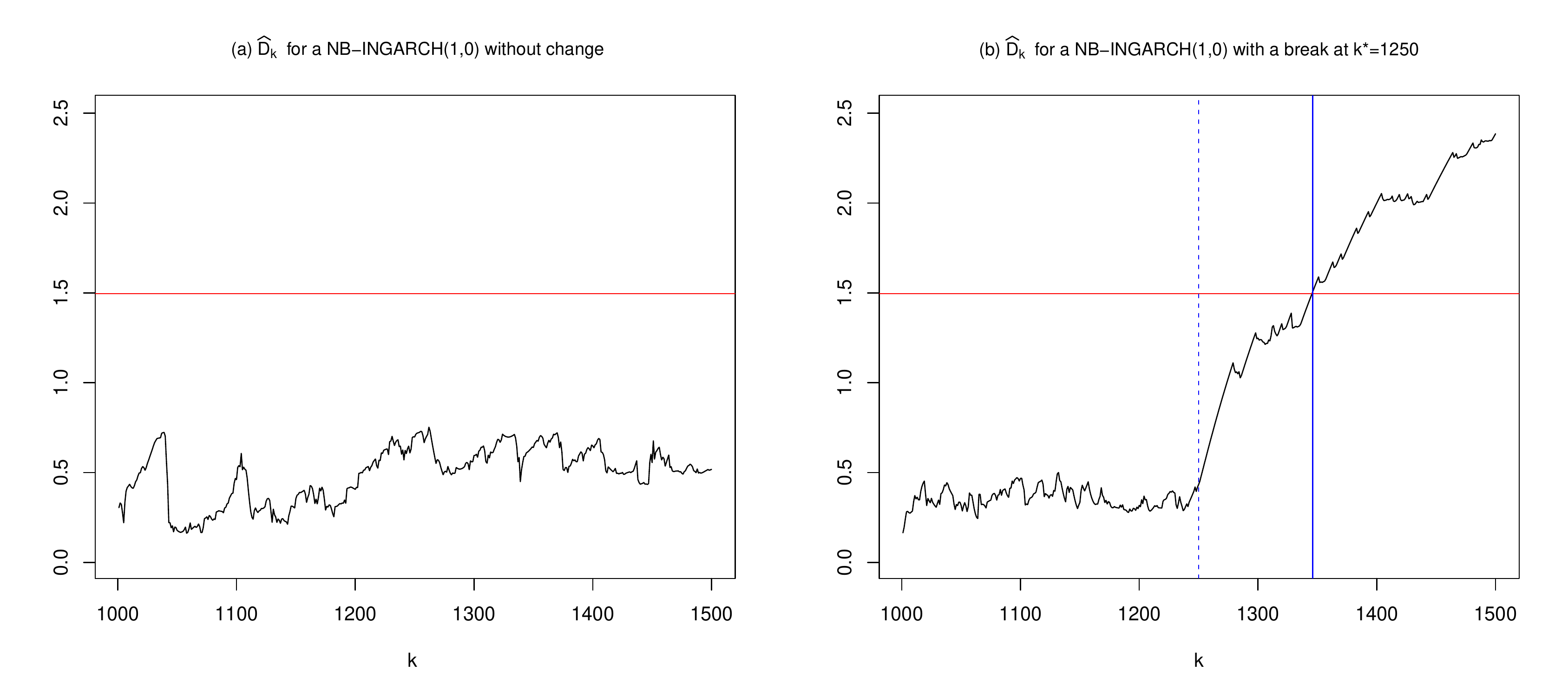}
\end{center}
\vspace{-.8cm}
\caption{\it \small Typical realization of the detector $\widehat{D}_{k}$ for the sequential change-point detection in two NB-INGARCH$(1,0)$ processes with $r=14$.
\textbf{(a)} the corresponding values for a NB-INGARCH$(1,0)$ process without change, where the parameter
$\theta^*_0=(0.5, 0.8)$ is constant. \textbf{(b)} the corresponding values for  a NB-INGARCH$(1,0)$ process where the parameter $\theta^*_0=(0.5, 0.8)$ changes
to $\theta^*_1=(0.15, 0.8)$ at $k^{*}=1250$.  
The horizontal line represents the limit of the critical region of the test.}
\label{Graphe_Online_NB}
\end{figure}

\medskip
For $m=1000$, Figure \ref{Graphe_Online_NB} displays a typical realization of the statistic $(\widehat{D}_{k})_{1001\leq k \leq 1500}$ from the model (\ref{NB_INGARCH})  in a scenario without change and a scenario with a change-point at $k^* =1250$.
As can be seen in Figure \ref{Graphe_Online_NB}(a), in the scenario without change, the detector $\widehat{D}_{k}$ is under the horizontal line which represents the limit of the critical region. 
Figure \ref{Graphe_Online_NB}(b) shows that before change occurs, $\widehat{D}_{k}$ is less than the the limit of the critical region. But, after the break, the detector increases with a high speed until exceed the critical value;  such growth over a long period indicates that something (structural change) is happening in the model. 

Now, we consider scenarios under H$^*_0$ and H$^*_1$ with break at $k^* = 1.25m$ in the  model (\ref{NB_INGARCH}). 
For $r=3,14$ and $m = 150,500,1000$, Table \ref{Table.2} indicates the empirical levels and the empirical powers based on $200$ replications. The case where $m=150$ is related and close to the real data example. 
Some elementary statistics of the empirical detection delays are summarized in Table \ref{Table.3}. \\
  The results of Table \ref{Table.2} show some distortion in terms of the empirical level for moderate sample sizes (see for instance, $m=150$). However, one can see that the empirical levels of the procedure approaching the nominal level when $m$ increases. In addition, for all scenarios considered, the empirical powers increases with $m$ and tends to approach one as $m$ increases;  which is in accordance with the asymptotic results of Theorem \ref{th4} and Corollary \ref{corol2}.\\ 
  %
  In Table \ref{Table.3}, let us recall that the detection delay $\widehat{d}_{m}$ (defined in (\ref{delay})) is the random distance between the break instant and the stopping time of the procedure. For example, when $m= 150$ with a change-point occurred at the time $k^* = 187$; from Table\ref{Table.3}, this change-point is detected on average after a delay of 20, 23, 23 and 20 respectively for these scenarios.
 Also, one can see that, for two historical sample sizes $m_1$ and $m_2$ with $m_1 < m_2$, the sequence $\widehat{d}_{m_2} - \sqrt{m_2/m_1} \widehat{d}_{m_1} $ decreases when $m_1$ and $m_2$ increases. 
This is in accordance with Theorem \ref{th4} where $\widehat{d}_m$ can be bounded by  $\mathcal{O}_P\Big ( \min\big(m^{1/2 + \epsilon}, \sqrt{m}(\log m)^a \big) \Big)$ for any $\epsilon, a>0$. 

\begin{table}[h!]
\centering
\scriptsize
\footnotesize
\caption{\small \it Empirical levels and powers at the nominal level $0.05$ for the sequential change-point detection in a NB-INGARCH process.}
\label{Table.2}
\vspace{.2cm}
\hspace*{-.3cm}
\begin{tabular}{lllccc}
\hline
\rule[0cm]{0cm}{.5cm}
&&$r$&$m=150$&$m=500$&$m=1000$\\
\hline
\rule[0cm]{0cm}{.4cm}
 Empirical levels: &&&&&\\

                  &$\theta^*_0=(0.5, 0.8)$   &$3$  &$0.095$&$0.035$&$0.040$\\
                  &                          &$14$ &$0.120$&$0.040$&$ 0.055$\\

& & &&&\\
                  & $\theta^*_0=(8.2, 0.2, 0.13)$ &$3$  &$0.140$&$0.080$&$0.065$\\
                  &                               &$14$ &$0.150$&$0.070$&$0.060$\\

&&&&&\\
\rule[0cm]{0cm}{.2cm}
Empirical powers: & & &&\\
                   &$\theta^*_0=(0.5, 0.8);~~~\theta^*_1=(0.15, 0.8);$  &$3$  &$0.305$&$0.625$&$0.985$\\
                  &                                                     &$14$ &$0.335$&$0.740$&$0.990$\\

& &&&&\\

                  & $\theta^*_0=(8.2, 0.2, 0.13); ~~~\theta^*_1=(2.35, 0.12, 0.61);$ &$3$  &$0.420$&$0.545$&$0.755$\\
                  &                                                                  &$14$ &$0.435$&$0.640$&$0.880$\\

 \Xhline{.75pt}
\end{tabular}
\end{table}

\begin{table}[h!]
\scriptsize
\footnotesize
\centering
\caption{ \it \small Some elementary statistics of the empirical detection delay for sequential change-point detection in a NB-INGARCH(1,0) process with $\theta^*_0=(0.5, 0.8),~\theta^*_1=(0.15, 0.8)$, and a NB-INGARCH(1,1) process with $\theta^*_0=(8.2, 0.2, 0.13), ~\theta^*_1=(2.35, 0.12, 0.61)$.}
\label{Table.3}
\vspace{.2cm}
\begin{tabular}{lllcccccccc}
\hline
 \rule[0cm]{0cm}{.35cm}
  &&&&\multicolumn{7} {l} {$\widehat{d}_{m}$}  \\
\cline{5-11}
\rule[0cm]{0cm}{.3cm}
         &&&&Mean&SD&Min&$Q_1$&Med&$Q_3$&Max\\
\Xhline{.6pt}
 \rule[0cm]{0cm}{.35cm}
NB-INGARCH$(1,0)$:
                 &$r=3$&$m=150;~~k^*=187$    &&$19.80$&$9.40$&$0$&$14$&$19$&$26$&$35$\\
                 \rule[0cm]{0cm}{.25cm}
                 &&$m=500;~~k^*=625$        &&$94.05$&$21.55$&$35$&$79$&$98$&$113$&$125$\\ 
                 \rule[0cm]{0cm}{.25cm}
                 &&$m=1000;~k^*=1250$      &&$142.7$&$46.03$&$19$&$115$&$141$&$176$&$246$\\
                 
                 \rule[0cm]{0cm}{.5cm}
                 &$r=14$&$m=150;~~k^*=187$  &&$22.57$&$10.98$&$3$&$11$&$27$&$31$&$38$\\
                 \rule[0cm]{0cm}{.25cm}
                 &&$m=500;~~k^*=625$        &&$84.78$&$23.09$&$10$&$70$&$88$&$103$&$123$\\
                 \rule[0cm]{0cm}{.25cm}
                 &&$m=1000;~k^*=1250$      &&$129.9 $&$42.40$&$25$&$98$&$126$&$157$&$234$\\  
                  
\rule[0cm]{0cm}{.8cm}

NB-INGARCH$(1,1)$: 
                &$r=3$&$m=150;~~k^*=187$     &&$23.31$&$11.42$&$1$&$18$&$27$&$33$&$38$\\
                 \rule[0cm]{0cm}{.25cm}
                 &&$m=500;~~k^*=625$         &&$78.55$&$32.18$&$8$&$57$&$83$&$103$&$125$\\
                 \rule[0cm]{0cm}{.25cm}
                 &&$m=1000;~k^*=1250$       &&$145.1$&$64.71$&$3$&$92$&$157$&$197$&$244$\\ 
                 
                 \rule[0cm]{0cm}{.5cm}
                 &$r=14$&$m=150;~~k^*=187$    &&$20.17$&$11.54$&$0$&$11$&$19$&$31$&$38$\\
                 \rule[0cm]{0cm}{.25cm}
                 &&$m=500;~~k^*=625$          &&$85.14$&$34.69$&$4$&$59$&$97$&$113$&$125$\\
                 \rule[0cm]{0cm}{.25cm}
                 &&$m=1000;~k^*=1250$        &&$129.2$&$57.54$&$20$&$88$&$127$&$166$&$245$\\

 \Xhline{.9pt}
\end{tabular}
\end{table}


\subsection{Real data application}
\begin{figure}[h!]
\begin{center}
\includegraphics[height=6cm, width=17cm]{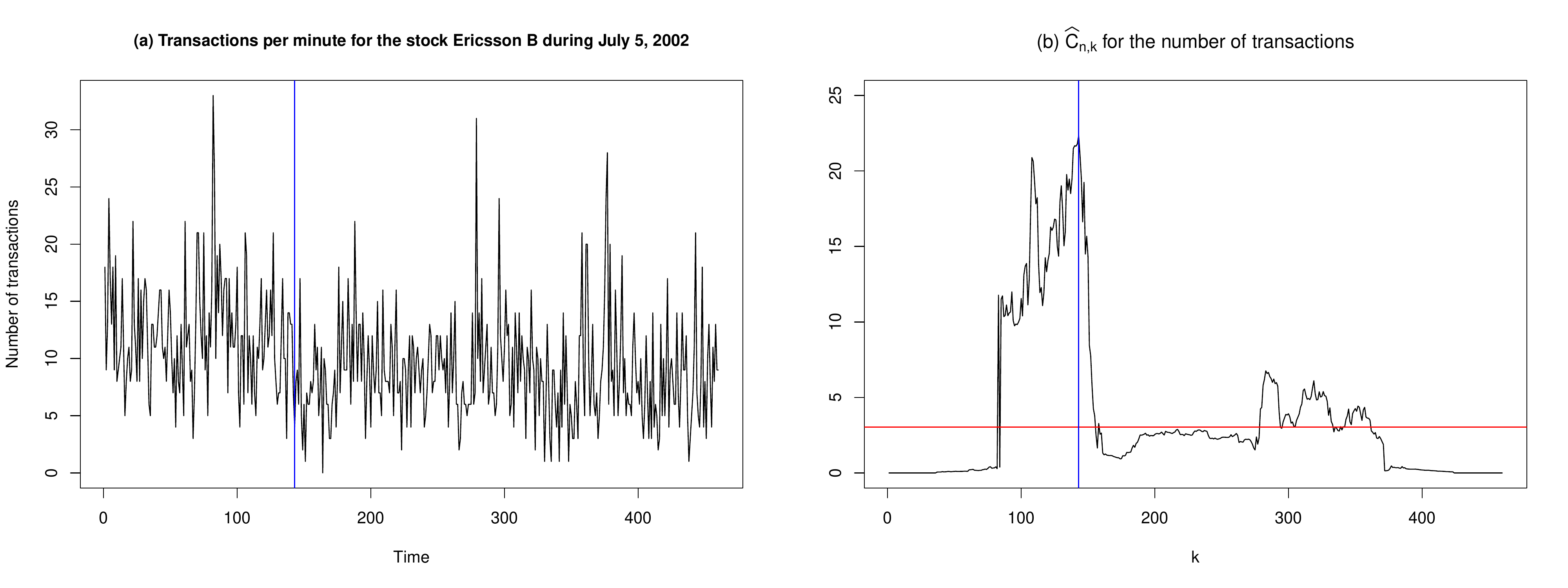} 
\end{center}
\vspace{-.8cm}
\caption{\small \it Plot of the statistic $\widehat{C}_{n,k}$ for the retrospective change-point detection applied to the number of transactions per minute for the stock Ericsson B during July 5, 2002 with an INGARCH(1,1) representation. 
The horizontal line in \textbf{(b)} represents the limit of the critical region of the test and the vertical line is the estimated breakpoint.}
\label{Graphe_Offline_Ericson}
\end{figure}

Consider the series of the number of transactions per minute for the stock Ericsson B during July 5, 2002 (see Figure \ref{Graphe_Offline_Ericson}(a)). 
 There are $460$ available observations that represent the transaction from $09:35$ through $17:14$. 
The empirical mean is $9.824$ while the empirical variance is $23.753$, which indicates that the data are overdispersed. 
See for instance, Fokianos \textit{et al.} (2009), Fokianos and Neumann (2013), Davis and Liu (2016),  Doukhan and Kengne (2015), Diop and Kengne (2017)  for some works that carried out an application to such data. 
 This series has been already analyzed by Diop and Kengne (2017) with a retrospective change point procedure based on the  maximum likelihood estimator of the model's parameter, under the assumption that, the conditional distribution of the data is negative binomial. We carry out an application without this assumption.  
 
 \medskip
Firstly, we apply the off-line change-point detection procedure with a INGARCH$(1,1)$ representation. 
 The realizations of the statistic $\widehat{C}_{n,k}$ with $q \equiv1 $, $u_n=[\left(\log(n)\right)^{2.5}]$ and $v_n=[\left(\log(n)\right)^{2}]$ are displayed in Figure \ref{Graphe_Offline_Ericson}(b). A change point is found at $\widehat{t}=143$, which exactly corresponds to the break that has been detected by Diop and Kengne (2017) (under the negative binomial assumption). 
 With the PQMLE, the estimated model on each regime yields: 
\begin{equation*}
\widehat \lambda_{t}=\left\{
\begin{array}{ll}
  \underset{(2.69)}{8.35} + \underset{(0.04)}{0.20}  Y_{t-1} + \underset{(0.23)}{0.11} \widehat \lambda_{t-1}~\text{ for }~t\leq 143,  
  \\
\\
\underset{(0.63)}{2.35} + \underset{(0.02)}{0.12}  Y_{t-1} + \underset{(0.08)}{0.61}  \widehat \lambda_{t-1}~\text{ for }~t > 143  , \\
\end{array}
\right.
\end{equation*}
where in parentheses are the standard errors of the estimators obtained from the robust sandwich matrix.  

\begin{figure}[h!]
\begin{center}
\includegraphics[height=6.5cm, width=13cm]{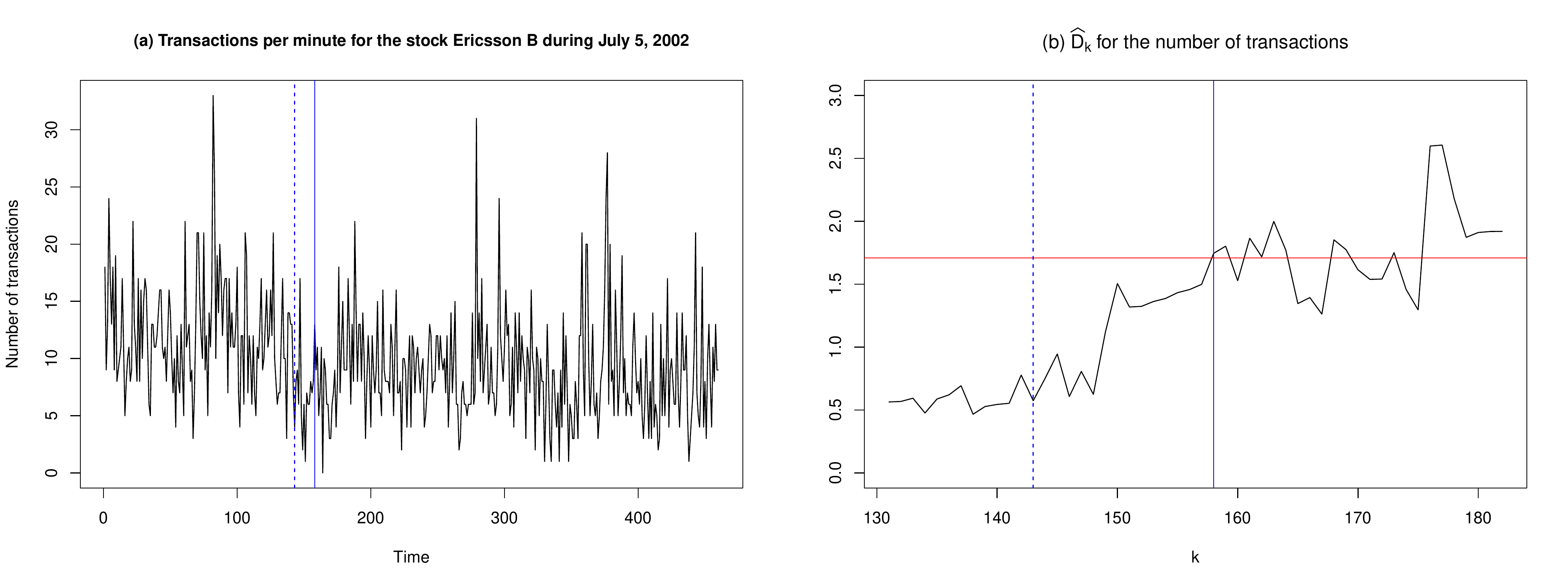} 
\end{center}
\vspace{-.8cm}
\caption{ \small \it Plot of the detector $\widehat{D}_{k}$ for the sequential change-point detection applied to  the number of transactions per minute for the stock Ericsson B during July 5, 2002 with an INGARCH(1,1) representation.
 The horizontal line in \textbf{(b)} represents the limit of the critical region of the test, 
the dotted  line represents the break that has been detected from the retrospective procedure and 
the solid line indicates the stopping time of the sequential procedure.}
\label{Graphe_Online_Ericson}
\end{figure}

Secondly, we apply the sequential procedure to this series by considering the closed-end setting with $T=1.5$ and the observations from $t=1$ to $t=130$ as the historical data. Therefore, the monitoring starts at the time $t=131$ and may continue until the time $t=195$ if no break is detected before this instant; i.e, the monitoring period is $\{131, 132,\cdots,195\}$. 
Figure \ref{Graphe_Online_Ericson} displays the realizations of the detector $\widehat{D}_{k}$ for $k=131,\cdots,182$.
 From this figure, one can see that the sequential procedure stops at time $k=158$; that is, $ \widehat d_m=158-143=15$ minutes after the break time detected from the retrospective procedure, which is reasonably good.
 
\section{Concluding remarks}
 This contribution addresses together the retrospective and the sequential change-point detection in a general class of integer-valued time series.
Numerous works have been done on these directions by assuming that the conditional distribution of the process is known; which is quite restrictive for practical issues. 
To overcome this drawback, we tackle these questions in a semiparametric framework with procedures based on the Poisson QMLE. 
For both the retrospective and the sequential detection, we propose test statistics that converge to some distributions obtained from the standard Brownian motion under the null hypothesis of no change and
diverge to infinity under the alternative; that is, these procedures are consistent. 
In the sequential detection,  the updated estimator which is computed without the historical observations leads to a procedure with a reasonably good detection delay that can be  bounded by $\mathcal{O}_P\Big ( \min\big(m^{1/2 + \epsilon}, \sqrt{m}(\log m)^a \big) \Big)$ for any $\epsilon, a>0$.  
Empirical studies show that these procedures overall work well for simulated and real data example. 
 A good extension of this work is to carry out these procedures with the estimator that are based on the negative binomial QMLE (see Aknouche {\it et al.} (2018)).


 \section{Proofs of the main results}  
  Let $ (\psi_n)_{n \in \N}  $ and $ (r_n)_{n \in \N}  $ be sequences of random variables or vectors. Throughout this section, we use the notation
  $ \psi_n = o_P(r_n)  $ to mean:  for all $  \varepsilon  > 0, ~ \prob( \|\psi_n \| \geq \varepsilon \|r_n \| ) \limiten 0$.
  Write $ \psi_n = O_P(r_n)$ to mean:  for all  $  \varepsilon > 0 $,  there exists $C>0$  such that  $\prob( \|\psi_n \| \geq C \|r_n \| )\leq \varepsilon $
   for $n$ large enough. 
   In the sequel, $C$ denotes a positive constant  whose the value may differ from one inequality to
another.

 \subsection{Results of the retrospective change-point detection} 
 The following lemma is obtained from the Lemma A.1 and A.4 of Diop et Kengne (2021) (for (i.) see also the proof of Theorem 2.1 of Ahmad and Francq (2016)); the proof is then omitted.  
 \begin{lem}\label{Lem_0} 
 Assume that the assumptions of Theorem \ref{th1} hold. 
 %
 %
Then, 
  \[
(i.)~\frac{1}{n}\left\| \widehat{L}(T_{1,n},\theta)-L(T_{1,n},\theta)\right\|_{\Theta} \limiten 0
	~\text{ and }~
(ii.)~ \frac{1}{\sqrt{n}}\Big\| \frac{\partial \widehat{L}(T_{1,n},\theta)}{\partial \theta}  -\frac{\partial L(T_{1,n},\theta)}{\partial \theta}\Big\|_{\Theta} \limiten 0.
\]
 \end{lem} 
 %
 
%
 \subsubsection{Proof of Theorem \ref{th1}}
  Define the statistic
$$
C_n=\max_{v_n<k<n-v_n}C_{n,k},~~~~~~~~~~~~~~~~~~~~~~~~~~~~~~~~~~~~~~~~~~~~~~~~~~~~~~~~~~$$
with
 \[ C_{n,k}=\frac{1}{q^{2}(\frac{k}{n})}\frac{k^{2}(n-k)^{2}}{n^{3}}\left(\widehat{\theta}(T_{1,k})-\widehat{\theta}(T_{k+1,n})\right)' \Sigma \left(\widehat{\theta}(T_{1,k})-\widehat{\theta}(T_{k+1,n})\right);
\]
where $\Sigma$ is defined at (\ref{res_Francq}) and computed at $\theta^*_1$.
Let $k,k^{\prime} \in[1,n]$, $\bar{\theta} \in \Theta$ and $i \in \{1,2,\cdots,d\}$. The Taylor expansion to the function $\theta \mapsto \frac{\partial}{\partial \theta_i} L({T_{k,k^{\prime}}},\theta)$ implies that there exists $\theta_{n,i}$ between $\bar{\theta}$ and $\theta^*_{1}$ such that
\[ 
\frac{\partial}{\partial \theta_i} L({T_{k,k^{\prime}}},\bar{\theta})=\frac{\partial}{\partial \theta_i} L({T_{k,k^{\prime}}},\theta_0) +\frac{\partial^{2}}{\partial \theta\partial \theta_i} L({T_{k,k^{\prime}}},\theta_{n,i})(\bar{\theta}-\theta_0).
\]
It is equivalent to
\begin{eqnarray}\label{Eq_Taylor}
(k^{\prime}-k+1)J_n({T_{k,k^{\prime}}},\bar{\theta}).(\bar{\theta}-\theta_0)=\frac{\partial}{\partial \theta} L({T_{k,k^{\prime}}},\theta_0)-\frac{\partial}{\partial \theta} L({T_{k,k^{\prime}}},\bar{\theta}),
\end{eqnarray}
where
\begin{eqnarray}\label{G_n}
J_n({T_{k,k^{\prime}}},\bar{\theta})= - \frac{1}{(k^{\prime}-k+1)}\frac{\partial^{2}}{\partial \theta\partial \theta_i} L({T_{k,k^{\prime}}},\theta_{n,i})_{1\leq i \leq d}.
\end{eqnarray}
%
 %
 The following lemma will be useful in the sequel.
\begin{lem}\label{Lem_1}
Suppose that the assumptions of Theorem \ref{th1} hold. Then,
\begin{enumerate}
    \item  $\max_{v_n<k<n-v_n}\big|\widehat{C}_{n,k}-C_{n,k}\big|=o_P(1)$;
    
    \item $\left(\frac{\partial}{\partial \theta} \ell_t(\theta^*_1),\mathcal{F}_{t}\right)_{t \in \mathbb{Z}}$ is a stationary ergodic, square integrable martingale difference sequence with covariance matrix $I$ (computed at $\theta^*_1$);
    
\item $\mathbb{E}\left(\frac{\partial^{2} \ell_0 (\theta^*_1)}{\partial \theta\partial \theta'}\right)=-J$ (computed at $\theta^*_1$);

\item $J_n({T_{1,n}},\widehat{\theta}({T_{1,n}})) \limitepsn J$.
\end{enumerate}
\end{lem}

\noindent
{\bf Proof.} 
\begin{enumerate}
\item See the proofs of Lemma {\bf A.1} of Diop and Kengne (2017) and Lemma 7.3. of Doukhan and Kengne (2015); by using the same arguments, one can go along similar lines as in these proofs.

\item 
Under H$_0$, $(X_t,Y_t)_{t \in \Z}$ is a stationary and ergodic process, the same properties hold for
 $\left(\frac{\partial}{\partial \theta} \ell_t(\theta^*_1)\right)_{t \in \mathbb{Z}}$.
 Moreover, 
\begin{equation}\label{def.1_ellt}
\frac{\partial \ell_t(\theta)}{\partial\theta}=\Big(\frac{Y_{t}}{\lambda_t(\theta)}-1 \Big) \frac{\partial \lambda_t(\theta)}{\partial\theta},~ \text{ for all } \theta \in \Theta.
\end{equation}
Since $\lambda_t(\theta)$ and $ \frac{\partial \lambda_t(\theta)}{\partial\theta}$ are $\mathcal{F}_{t-1}$-measurable for any $\theta \in \Theta$, it holds that
$\mathbb{E}\left(\frac{\partial \ell_t(\theta^*_1)}{\partial\theta}|\mathcal{F}_{t-1}\right)=0$.\\
Also,
\[
\mathbb{E}\left[ \left(\frac{\partial \ell_0(\theta^*_1)}{\partial\theta}\right)\left(\frac{\partial \ell_0(\theta^*_1)}{\partial\theta}\right)' \right]= \mathbb{E}\left[ \Big(\frac{Y_{0}}{\lambda_t(\theta^*_1)}-1 \Big)^2 \frac{\partial \lambda_t(\theta^*_1)}{\partial\theta} \frac{\partial \lambda_t(\theta^*_1)}{\partial\theta'}\right]=I.
\]

\item 
From (\ref{def.1_ellt}), we deduce 
\begin{equation}\label{def.2_ellt}
\frac{\partial^{2}\ell_t(\theta)}{\partial \theta\partial\theta'}
= \left(\frac{Y_{t}}{\lambda_t(\theta)}-1 \right)\frac{\partial^{2}\lambda_t(\theta)}{\partial \theta\partial\theta'} -
\frac{Y_t}{\lambda_t(\theta)^2} \frac{\partial \lambda_t(\theta)}{\partial\theta} \frac{\partial \lambda_t(\theta)}{\partial\theta'},~ \text{ for all } \theta \in \Theta.
\end{equation}
Since $\frac{\partial^{2}\lambda_t(\theta)}{\partial \theta\partial\theta'}$ is also  $\mathcal{F}_{t-1}$-measurable for any $\theta \in \Theta$,
 it holds that $\E\Big[\left(\frac{Y_{t}}{\lambda_t(\theta^*_1)}-1 \right)\frac{\partial^{2}\lambda_t(\theta^*_1)}{\partial \theta\partial\theta'} | \mathcal{F}_{t-1}\Big]=0$.\\ 
Hence, 
\[\E\Big(\frac{\partial^{2}\ell_t(\theta^*_1)}{\partial \theta\partial\theta'}\Big)=-\E\left(\frac{Y_t}{\lambda_t(\theta^*_1)^2} \frac{\partial \lambda_t(\theta^*_1)}{\partial\theta} \frac{\partial \lambda_t(\theta^*_1)}{\partial\theta'}\right)
= -\E\left(\frac{1}{\lambda_t(\theta^*_1)} \frac{\partial \lambda_t(\theta^*_1)}{\partial\theta} \frac{\partial \lambda_t(\theta^*_1)}{\partial\theta'}\right) =-J.\]

\item By applying $(\ref{G_n})$ with $\bar{\theta} = \widehat{\theta}(T_{1,n})$, it holds that
\begin{eqnarray*}
J_n(T_{1,n},\widehat{\theta}(T_{1,n}))= - \Big(\frac{1}{n}\frac{\partial^{2}}{\partial \theta\partial \theta_i} L(T_{1,n},\theta_{n,i})\Big)_{1\leq i \leq d}=-\frac{1}{n}\Big(\sum_{t=1}^{n}\frac{\partial^{2}}{\partial \theta\partial \theta_i} \ell_t(\theta_{n,i})\Big)_{1\leq i \leq d},
\end{eqnarray*}
where $\theta_{n,i}$ belongs between $\widehat{\theta}(T_{1,n})$ and $\theta_1^*$. 
 Since $\widehat{\theta}(T_{1,n}), \theta_{n,i} ~ (\text{for any } i=1,\cdots,d)  \limitepsn \theta_1^*$, from the proof of Theorem 2.2 of \cite{Francq2016}, we get
\[ -\frac{1}{n}\Big(\sum_{t=1}^{n}\frac{\partial^{2}}{\partial \theta\partial \theta_i} \ell_t(\theta_{n,i})\Big)_{1\leq i \leq d} \limitepsn  -\E\Big(\frac{\partial^{2}\ell_0(\theta^*_1)}{\partial \theta\partial\theta'}\Big) = J. \]

\end{enumerate}
\begin{flushright}
$\Box$ 
\end{flushright}

Now, let us use Lemma \ref{Lem_1} to show that
\[C_n\stackrel{\mathcal{D}}{\longrightarrow}\sup_{0\leq\tau\leq1}\frac{\left\|W_d(\tau)\right\|^{2}}{q^{2}(\tau)}~~ \textrm{as}~~n\rightarrow \infty.\]
\noindent Let $v_n \leq k \leq n-v_n$. By applying (\ref{Eq_Taylor}) with $\bar{\theta}=\widehat{\theta}(T_{1,k})$ and $T_{k,k^\prime}=T_{1,k}$, we get
\begin{eqnarray}\label{Eq_T_1,k}
J_n(T_{1,k},\widehat{\theta}(T_{1,k})) \cdot (\widehat{\theta}(T_{1,k})-\theta_1^*)=\frac{1}{k}\left(\frac{\partial}{\partial \theta} L(T_{1,k},\theta_1^*)-\frac{\partial}{\partial \theta} L(T_{1,k},\widehat{\theta}(T_{1,k}))\right).
\end{eqnarray}
With $\bar{\theta}=\widehat{\theta}(T_{k+1,n})$ and $T_{k,k^\prime}=T_{k+1,n}$, (\ref{Eq_Taylor}) becomes
\begin{eqnarray}\label{Eq_T_k+1,n}
J_n(T_{k+1,n},\widehat{\theta}(T_{k+1,n})) \cdot (\widehat{\theta}(T_{k+1,n})-\theta_1^*)=\frac{1}{n-k}\left(\frac{\partial}{\partial \theta} L(T_{k+1,n},\theta_1^*)-\frac{\partial}{\partial \theta} L(T_{k+1,n},\widehat{\theta}(T_{k+1,n}))\right).
\end{eqnarray}
As $n\rightarrow +\infty$, we have
\begin{align*}
&\left\|J_n(T_{1,k},\widehat{\theta}(T_{1,k}))-J\right\|=o(1),~~~\left\|J_n(T_{k+1,n},\widehat{\theta}(T_{k+1,n}))-J\right\|=o(1)\\ 
&
\sqrt{k}\left(\widehat{\theta}(T_{1,k})-\theta_1^*\right)=O_P(1) ~~\text{and}~~ \sqrt{n-k}\left(\widehat{\theta}(T_{k+1,n})-\theta_1^*\right)=O_P(1).
\end{align*}
According to (\ref{Eq_T_1,k}), for  $n$ large enough, we get
\begin{align*}
\sqrt{k}J\left(\widehat{\theta}(T_{1,k})-\theta_1^*\right)
&=
\frac{1}{\sqrt{k}}\left(\frac{\partial}{\partial \theta} L(T_{1,k},\theta_1^*)-\frac{\partial}{\partial \theta} L(T_{1,k},\widehat{\theta}(T_{1,k}))\right) 
\\
& \hspace{6.5cm} -\sqrt{k}\left(\left(J_n(T_{1,k},\widehat{\theta}(T_{1,k}))-J\right)\left(\widehat{\theta}(T_{1,k})-\theta_0\right)\right)\\
&=\frac{1}{\sqrt{k}}\left(\frac{\partial}{\partial \theta} L(T_{1,k},\theta_1^*)-\frac{\partial}{\partial \theta} L(T_{1,k},\widehat{\theta}(T_{1,k}))\right) +o_P(1)\\
&= \frac{1}{\sqrt{k}}\left(\frac{\partial}{\partial \theta} L(T_{1,k},\theta_1^*)- \frac{\partial}{\partial \theta} \widehat L(T_{1,k},\widehat{\theta}(T_{1,k}))\right)+o_P(1)\\
& \hspace{6cm} +\frac{1}{\sqrt{k}}\left(\frac{\partial}{\partial \theta} \widehat L(T_{1,k},\widehat{\theta}(T_{1,k}))-\frac{\partial}{\partial \theta} L(T_{1,k},\widehat{\theta}(T_{1,k}))\right) \\
&=\frac{1}{\sqrt{k}}\left(\frac{\partial}{\partial \theta} L(T_{1,k},\theta_1^*)-\frac{\partial}{\partial \theta} \widehat L(T_{1,k},\widehat{\theta}(T_{1,k}))\right) +o_P(1)~ ~(\text{from Lemma}~ \ref{Lem_0}) .
\end{align*}
It is equivalent to
\begin{equation}\label{Eq_a}
J\left(\widehat{\theta}(T_{1,k})-\theta_1^* \right)=\frac{1}{k}\left(\frac{\partial}{\partial \theta} L(T_{1,k},\theta_1^*)-\frac{\partial}{\partial \theta} \widehat L(T_{1,k},\widehat{\theta}(T_{1,k}))\right) +o_P\left(\frac{1}{\sqrt{k}}\right).
\end{equation}
 For $n$ large enough, $\widehat{\theta}(T_{1,k})$ is an interior point of $\Theta$ and we have $ \frac{\partial}{\partial \theta} \widehat L(T_{1,k},\widehat{\theta}(T_{1,k}))=0$.
Hence, for $n$ large enough, we get from (\ref{Eq_a})
\begin{equation}\label{Eq_a_bis}
J\left(\widehat{\theta}(T_{1,k})-\theta_1^*\right)
 =\frac{1}{k}\frac{\partial}{\partial \theta} L(T_{1,k},\theta_1^*)+o_P\left(\frac{1}{\sqrt{k}}\right).
\end{equation}
Similarly, we can use (\ref{Eq_T_k+1,n}) to obtain
\begin{equation}\label{Eq_a_bisbis}
J\left(\widehat{\theta}(T_{k+1,n})-\theta_1^*\right)=\frac{1}{n-k}\frac{\partial}{\partial \theta} L(T_{k+1,n},\theta_1^*)+o_P\left(\frac{1}{\sqrt{n-k}}\right).
\end{equation}
The subtraction of the two above equalities gives
\begin{align*}
J\left(\widehat{\theta}(T_{1,k})-\widehat{\theta}(T_{k+1,n})\right)&=\frac{1}{k}\frac{\partial}{\partial \theta} L(T_{1,k},\theta_1^*)-\frac{1}{n-k}\frac{\partial}{\partial \theta} L(T_{k+1,n},\theta_1^*)+o_P\left(\frac{1}{\sqrt{k}}+\frac{1}{\sqrt{n-k}}\right)\\
&=\frac{1}{k}\frac{\partial}{\partial \theta} L(T_{1,k},\theta_1^*)-\frac{1}{n-k}\left(\frac{\partial}{\partial \theta} L(T_{1,n},\theta_1^*)-\frac{\partial}{\partial \theta} L(T_{1,k},\theta_1^*)\right)
 +o_P\left(\frac{1}{\sqrt{k}}+\frac{1}{\sqrt{n-k}}\right)\\
&=\frac{n}{k(n-k)}\left(\frac{\partial}{\partial \theta} L(T_{1,k},\theta_1^*)-\frac{k}{n}.\frac{\partial}{\partial \theta} L(T_{1,n},\theta_1^*)\right)+o_P\left(\frac{1}{\sqrt{k}}+\frac{1}{\sqrt{n-k}}\right);
\end{align*}
i.e.,
\begin{align*}
\frac{k(n-k)}{n^{3/2}}J(\widehat{\theta}(T_{1,k})-\widehat{\theta}(T_{k+1,n}))&=\frac{1}{\sqrt{n}}\left(\frac{\partial}{\partial \theta} L(T_{1,k},\theta_1^*)-\frac{k}{n}.\frac{\partial}{\partial \theta} L(T_{1,n},\theta_1^*)\right) +o_P\Bigg(\frac{\sqrt{k(n-k)}}{n} + \frac{ \sqrt{n-k}}{\sqrt{n}}\Bigg)\\
&=\frac{1}{\sqrt{n}}\left(\frac{\partial}{\partial \theta} L(T_{1,k},\theta_1^*)-\frac{k}{n}.\frac{\partial}{\partial \theta} L(T_{1,n},\theta_1^*)\right) +o_P(1).
\end{align*}
According the above equation, we have
\begin{equation}\label{Eq_b}
\frac{k(n-k)}{n^{3/2}}I^{-1/2}J(\widehat{\theta}(T_{1,k})-\widehat{\theta}(T_{k+1,n}))=\frac{I^{-1/2}}{\sqrt{n}}\left(\frac{\partial}{\partial \theta} L(T_{1,k},\theta_1^*)-\frac{k}{n}.\frac{\partial}{\partial \theta} L(T_{1,n},\theta_1^*)\right) +o_P(1).
\end{equation}
Recall that for any $0<\tau<1$,  
\[
\frac{\partial}{\partial \theta}L(T_{1,[n \tau]},\theta_1^*)=\sum_{t=1}^{[n \tau]}\frac{\partial}{\partial \theta}\ell_t(\theta_1^*).
\]
%
The process $\left(\frac{\partial}{\partial \theta} \ell_t(\theta_1^*),\mathcal{F}_{t}\right)_{t \in \mathbb{Z}}$ is a stationary ergodic square integrable martingale difference process with covariance matrix $I$ (see Lemma \ref{Lem_1}). By applying the central limit theorem for the martingale difference sequence (see Billingsley (1968)), we have
\begin{align*}
\frac{1}{\sqrt{n}}\left(\frac{\partial}{\partial \theta}L(T_{1,[n \tau]},\theta_1^*)-\frac{[n \tau]}{n}\frac{\partial}{\partial \theta}L(T_{1,n },\theta_1^*)\right)
&=
\frac{1}{\sqrt{n}}\Bigg(\sum_{t=1}^{[n \tau]}\frac{\partial}{\partial \theta}\ell_t(\theta_1^*)-\frac{[n \tau]}{n}\sum_{t=1}^{n }\frac{\partial}{\partial \theta}\ell_t(\theta_1^*)\Bigg)\\
&~~
\limiteloin B_{I}(\tau)-\tau B_{I}(1),
\end{align*}
where $B_{I}$ is a Gaussian process with covariance matrix $\min(s,t)I$.\\
Hence,
\begin{align*}
\frac{1}{\sqrt{n}}I^{-1/2}\left(\frac{\partial}{\partial \theta}L(T_{1,[n \tau]},\theta_1^*)-\frac{[n \tau]}{n}\frac{\partial}{\partial \theta}L(T_{1,n },\theta_1^*)\right)~~
\limiteloin  B_{d}(\tau)-\tau B_{d}(1)=W_d(\tau) 
\text{ in } ~\mathcal{D}\left(\left[0,1\right]\right),
\end{align*}
 where $B_d$ is a $d$-dimensional standard motion, and $W_d$ is a $d$-dimensional Brownian bridge.\\
From (\ref{Eq_b}), as $n\rightarrow +\infty$, we have
\begin{align*}
C_{n,\left[n\tau\right]}&=
\frac{1}{q^{2}\left(\frac{[n\tau]}{n}\right)}\frac{[n\tau]^{2}(n-[n\tau])^{2}}{n^{3}}\left(\widehat{\theta}(T_{1,[n\tau]})-\widehat{\theta}(T_{[n\tau]+1,n})\right)'\Sigma \left(\widehat{\theta}(T_{1,[n\tau]})-\widehat{\theta}(T_{[n\tau]+1,n})\right)\\
&=\frac{1}{q^{2}\left(\frac{[n\tau]}{n}\right)}\Bigg\|\frac{[n\tau](n-[n\tau])}{n^{3/2}}I^{-1/2}J\left(\widehat{\theta}(T_{1,[n\tau]})-\widehat{\theta}(T_{[n\tau]+1,n})\right)\Bigg\|^{2}\\
&=\frac{1}{q^{2}\left(\frac{[n\tau]}{n}\right)}\Bigg\|\frac{I^{-1/2}}{\sqrt{n}}\left(\frac{\partial}{\partial \theta} L(T_{1,[n\tau]},\theta_1^*)-\frac{[n\tau]}{n}.\frac{\partial}{\partial \theta} L(T_{1,n},\theta_1^*)\right) +o_P(1)\Bigg\|^{2}\\
&=\frac{1}{q^{2}\left(\frac{[n\tau]}{n}\right)}\Bigg\|\frac{I^{-1/2}}{\sqrt{n}}\Big(\sum_{t=1}^{[n \tau]}\frac{\partial}{\partial \theta}\ell_t(\theta_1^*)-\frac{[n \tau]}{n}\sum_{t=1}^{n }\frac{\partial}{\partial \theta}\ell_t(\theta_1^*)\Big) \Bigg\|^{2} +o_P(1)\\
&~~~\stackrel{\mathcal{D}}{\longrightarrow} \frac{\left\|W_d(\tau)\right\|^{2}}{q^{2}(\tau)}~\textrm{in}~ \mathcal{D}\left(\left[0,1\right]\right).
\end{align*}
According to the properties of $q$, we have for any $0<\epsilon<1/2$
\begin{align*}
\max_{[n\epsilon]<k<n-[n\epsilon]}C_{n,k}&=\max_{[n\epsilon]<k<n-[n\epsilon]}\frac{1}{q^{2}\left(\frac{k}{n}\right)}\frac{k^{2}(n-k)^{2}}{n^{3}}\left(\widehat{\theta}(T_{1,k})-\widehat{\theta}(T_{k+1,n})\right)'\Sigma \left(\widehat{\theta}(T_{1,k})-\widehat{\theta}(T_{k+1,n})\right)\\
&=\sup_{\epsilon<\tau<1-\epsilon}\frac{1}{q^{2}\left(\frac{[n\tau]}{n}\right)}\frac{[n\tau]^{2}(n-[n\tau])^{2}}{n^{3}}\left(\widehat{\theta}(T_{1,[n\tau]})-\widehat{\theta}(T_{[n\tau]+1,n})\right)'\Sigma 
\left(\widehat{\theta}(T_{1,[n\tau]})-\widehat{\theta}(T_{[n\tau]+1,n})\right)\\
&=\sup_{\epsilon<\tau<1-\epsilon}\frac{1}{q^{2}\left(\frac{[n\tau]}{n}\right)}\Big\|\frac{[n\tau](n-[n\tau])}{n^{3/2}}
I^{-1/2}J
\left(\widehat{\theta}(T_{1,[n\tau]})-\widehat{\theta}(T_{[n\tau]+1,n})\right)\Big\|^{2}\\
&=\sup_{\epsilon<\tau<1-\epsilon}\frac{1}{q^{2}\left(\frac{[n\tau]}{n}\right)}\Big\|\frac{I^{-1/2}}{\sqrt{n}}\Big(\frac{\partial}{\partial \theta} L_n(T_{1,[n\tau]},\theta_1^*)-\frac{[n\tau]}{n}.\frac{\partial}{\partial \theta} L_n(T_{1,n},\theta_1^*)\Big) +o_P(1)\Big\|^{2}\\
&=\sup_{\epsilon<\tau<1-\epsilon}\frac{1}{q^{2}\left(\frac{[n\tau]}{n}\right)}\Big\|\frac{I^{-1/2}}{\sqrt{n}}\Big(\sum_{t=1}^{[n \tau]}\frac{\partial}{\partial \theta}\ell_t(\theta_1^*)-\frac{[n \tau]}{\sqrt{n}}\sum_{t=1}^{n }\frac{\partial}{\partial \theta}\ell_t(\theta_1^*)\Big) \Big\|^{2} +o_P(1)\\
&~~~\limiteloin \sup_{\epsilon<\tau<1-\epsilon}\frac{\left\|W_d(\tau)\right\|^{2}}{q^{2}(\tau)}. 
\end{align*}
Hence, we have shown that as $n\rightarrow +\infty$, 
\[
C_{n,\left[n\tau\right]}\stackrel{\mathcal{D}}{\longrightarrow} \frac{\left\|W_d(\tau)\right\|^{2}}{q^{2}(\tau)}~\textrm{in}~ \mathcal{D}\left(\left[0,1\right]\right)
\]
and for all $0<\epsilon<1/2$,
\[\max_{[n\epsilon]<k<n-[n\epsilon]}C_{n,k}=\sup_{\epsilon<\tau<1-\epsilon}C_{n,[n\tau]}\stackrel{\mathcal{D}}{\longrightarrow} \sup_{\epsilon<\tau<1-\epsilon}\frac{\left\|W_d(\tau)\right\|^{2}}{q^{2}(\tau)}.
\]
In addition, since $I(q,c)<+\infty$ for some $c>0$, one can show that (see also \cite{Csorgo_1986}) \[ \lim_{\tau\rightarrow 0}\frac{\left\|W_d(\tau)\right\|}{q(\tau)}<\infty~~\textrm{and}~~\lim_{\tau\rightarrow 1}\frac{\left\|W_d(\tau)\right\|}{q(\tau)}<\infty ~ a.s.\]
Hence, for $n$ large enough, we have
\[
C_n=\max_{v_n<k<n-v_n}C_{n,k}=\sup_{\frac{v_n}{n}<\tau<1-\frac{v_n}{n}}C_{n,[n\tau]} \limiteloin \sup_{0<\tau<1}\frac{\left\|W_d(\tau)\right\|^{2}}{q^{2}(\tau)}.\]
\begin{flushright}
$\blacksquare$ 
\end{flushright}

\subsubsection{Proof of Theorem \ref{th2}}
Assume that the trajectory $(Y_1,\cdots,Y_n)$ satisfies
\begin{equation} \label{Eq_H1}
Y_{t}=\left\{
\begin{array}{ll}
Y^{(1)}_{t}~~\textrm{for}~~t\leq t^{*},\\
\\
Y^{(2)}_{t}~~\textrm{for}~~t>t^{*},\\
\end{array}
\right.
\end{equation} \\
where $t^{*}=[\tau^* n]$ with $0<\tau^*<1$ and $\{Y^{(j)}_{t}, t \in \mathbb{Z}\}$ ($j=1,2$) is a stationary solution of the model (\ref{Model}) depending on $\theta^{*}_j$ with $\theta^{*}_1 \neq \theta^{*}_2$.\\
We have
\begin{align*}
 \widehat{C}_{n,t^{*}}
&=
\frac{1}{q^{2}(\frac{t^{*}}{n})}\frac{{t^{*}}^{2}(n-t^{*})^{2}}{n^{3}}\left(\widehat{\theta}(T_{1,t^{*}})-\widehat{\theta}(T_{t^{*}+1,n})\right)'\widehat{\Sigma}(u_n)\left(\widehat{\theta}(T_{1,t^{*}})-\widehat{\theta}(T_{t^{*}+1,n})\right) \\
\text{ and }~~~~~~~~&\\
\widehat{C}_{n}&=\max_{v_n\leq k \leq n-v_n}\widehat{C}_{n,k}\geq \widehat{C}_{n,t^{*}}.
\end{align*}
 Then, to prove the Theorem \ref{th2}, it suffices to show that $\widehat{C}_{n,t^{*}} \limiteproban +\infty$.\\
 %
Recall that the matrix used to construct the test statistic is $\widehat{\Sigma}(u_n)$ given by
\[
\widehat{\Sigma}(u_n)
=\frac{1}{2} 
\left[
\widehat J(T_{1,u_n})  \widehat I(T_{1,u_n})^{-1}   \widehat J(T_{1,u_n}) +
\widehat J(T_{u_n+1,n})  \widehat I(T_{u_n+1,n})^{-1}  \widehat J(T_{u_n+1,n})
\right].
\]
According to the asymptotic proprieties of the PQMLE, we have 
\[
\widehat{\theta}(T_{1,t^{*}}) \limitepsn \theta^{*}_1,~~~\widehat{\theta}(T_{1,u_n}) \limitepsn \theta^{*}_1~~~
\text{ and}~~~
\widehat J(T_{1,u_n})  \widehat I(T_{1,u_n})^{-1}   \widehat J(T_{1,u_n}) \limitepsn \Sigma^{(1)},
\]
where
\begin{equation*}
\Sigma^{(1)}= J_1  I^{-1}_1 J_1 ~\text{ with } ~
J_1 =\E \Big[ \frac{1}{\lambda_{0}(\theta^*_1 )}  \frac{\partial \lambda_{0}(\theta^*_1 )}{ \partial \theta} \frac{\partial \lambda_{0}(\theta^*_1 )}{ \partial \theta'}  \Big] 
~\text{ and } ~
        I_1 =\E \Big[ \Big(\frac{Y_{0}}{\lambda_{0}(\theta^*_1 )}-1\Big)^2 \frac{\partial \lambda_{0}(\theta^*_1 )}{ \partial \theta} \frac{\partial \lambda_{0}(\theta^*_1 )}{ \partial \theta'}  \Big]. 
\end{equation*}
%
%
 Moreover,  the asymptotic proprieties of the PQMLE implies 
$\widehat{\theta}(T_{t^{*}+1,n}) \limitepsn \theta^{*}_2$.
Recall that, by definition, the two matrices in the formula of $\widehat{\Sigma}_n(u_n)$ are positive semi-definite and the first one converges a.s. to $\Sigma^{(1)}$ which is positive definite. \\
Then, for $n$ large enough, we can write
\begin{align*}
\widehat{C}_{n}
&\geq \widehat{C}_{n,t^{*}}\\
&\geq \frac{1}{2}  \frac{1}{q^{2}(\frac{t^{*}}{n})}\frac{{t^{*}}^{2}(n-t^{*})^{2}}{n^{3}} \left(\widehat{\theta}(T_{1,t^{*}})-\widehat{\theta}(T_{t^{*}+1,n})\right)'\\
& \hspace{2cm} \Big[
\widehat J(T_{1,u_n})  \widehat I(T_{1,u_n})^{-1}   \widehat J(T_{1,u_n}) +
\widehat J(T_{u_n+1,n})  \widehat I(T_{u_n+1,n})^{-1}  \widehat J(T_{u_n+1,n})
\Big] \left(\widehat{\theta}(T_{1,t^{*}})-\widehat{\theta}(T_{t^{*}+1,n})\right)\\
&\geq  \frac{1}{2} \frac{1}{q^{2}(\frac{t^{*}}{n})}\frac{{t^{*}}^{2}(n-t^{*})^{2}}{n^{3}}\left(\widehat{\theta}(T_{1,t^{*}})-\widehat{\theta}(T_{t^{*}+1,n})\right)'
\Big[
\widehat J(T_{1,u_n})  \widehat I(T_{1,u_n})^{-1}   \widehat J(T_{1,u_n}) 
\Big]
\left(\widehat{\theta}(T_{1,t^{*}})-\widehat{\theta}(T_{t^{*}+1,n})\right)\\~\\
&\geq
\frac{1}{2} \frac{1}{  \underset{0<\tau \leq \tau^*}{\sup} q^{2}(\tau)}n\left(\tau^{*}(1-\tau^{*})\right)^{2}
\left(\widehat{\theta}(T_{1,t^{*}})-\widehat{\theta}(T_{t^{*}+1,n})\right)' \Big[
\widehat J(T_{1,u_n})  \widehat I(T_{1,u_n})^{-1}   \widehat J(T_{1,u_n}) 
\Big] \left(\widehat{\theta}(T_{1,t^{*}})-\widehat{\theta}(T_{t^{*}+1,n})\right)\\
&\geq
C \times n \left(\widehat{\theta}(T_{1,t^{*}})-\widehat{\theta}(T_{t^{*}+1,n})\right)'
\times
\left[
\widehat J(T_{1,u_n})  \widehat I(T_{1,u_n})^{-1}   \widehat J(T_{1,u_n}) 
\right]
\times
\left(\widehat{\theta}(T_{1,t^{*}})-\widehat{\theta}(T_{t^{*}+1,n})\right).
\end{align*}
Therefore, since 
\[
\widehat{\theta}(T_{1,t^{*}})-\widehat{\theta}(T_{t^{*}+1,n}) \limitepsn \theta^{*}_{1}-\theta^{*}_{2}\neq 0
~~~\text{and} ~~~
\widehat J(T_{1,u_n})  \widehat I(T_{1,u_n})^{-1}   \widehat J(T_{1,u_n}) \limitepsn \Sigma^{(1)};
\]
we deduce that, $\widehat{C}_{n}\limitepsn +\infty$.
This completes the proof of the theorem.
\begin{flushright}
$\blacksquare$
\end{flushright}

\subsection{Results of the sequential change-point detection}
%
For any $k > m$ and $\ell \in \Pi_{m,k}$, denote
\[
 D_{k,\ell}= \sqrt{m}\frac{k-\ell}{k} \big\| I^{-1/2} J(\widehat{\theta}(T_{\ell,k})-\widehat{\theta}(T_{1,m})) \big\|.
\]
where $I$ and $J$ are computed under H$^*_0$ and depend on $\theta^*_1$.
The following lemma will be useful.
\begin{lem}\label{Lem_2}
Under the assumptions of Theorem \ref{th3},
\begin{equation*}
\sup_{k>m}\max_{\ell \in \Pi_{m,k}} \frac{1}{b((k-\ell)/m)} \| \widehat{D}_{k,\ell}-D_{k,\ell}\|=o_P(1)
~\text{ as }~ m \rightarrow + \infty.
\end{equation*}
\end{lem} 

\noindent
{\bf Proof.}\\ 
Let $k > m$ and $\ell \in \Pi_{m,k}$. As $m \rightarrow + \infty$, from Ahmad and Francq (2016), it holds that 
$\| \widehat I(T_{1,m})^{-1/2} \widehat J(T_{1,m})-  I^{-1/2} J\| \stackrel{a.s}{\longrightarrow}0$, $ \|\widehat{\theta}(T_{1,m})-\theta^*_1\|=O_P(1/\sqrt{m})$  and for  $k>m$, $\|\widehat{\theta}(T_{\ell,k})-\theta^*_1\|=O_P(1/\sqrt{k-\ell+1})$.
Hence,
\begin{align*}
&\sup_{k>m}\max_{\ell \in \Pi_{m,k}} \frac{1}{b((k-\ell)/m)} \| \widehat{D}_{k,\ell}-D_{k,\ell}\|\\
&~~~~\leq
\frac{1}{\inf_{s>0}b(s)}
\sup_{k>m}\max_{\ell \in \Pi_{m,k}} \sqrt{m} \frac{k-\ell}{k}\big\| \big(\widehat I(T_{1,m})^{-1/2} \widehat J(T_{1,m})-I^{-1/2} J\big)\big(\widehat{\theta}(T_{\ell,k})-\widehat{\theta}(T_{1,m})\big)\big\|\\
&~~~~\leq
C \sqrt{m} \big\| \widehat I(T_{1,m})^{-1/2} \widehat J(T_{1,m})-I^{-1/2} J\big\| \big\|\widehat{\theta}(T_{1,m})-\theta^*_1\big\|\\
& \hspace{3cm} + C \sup_{k>m}\max_{\ell \in \Pi_{m,k}} \sqrt{\frac{m}{k}} \cdot
 \frac{k-\ell}{\sqrt{k}} \cdot \big\| \widehat I(T_{1,m})^{-1/2} \widehat J(T_{1,m})-I^{-1/2} J\big\| \big\|\widehat{\theta}(T_{\ell,k})-\theta^*_1\big\|
\\
&~~~~\leq o_P(1)+C \sup_{k>m}\max_{\ell \in \Pi_{m,k}}
 \sqrt{k-\ell+1}\cdot \big\| \widehat I(T_{1,m})^{-1/2} \widehat J(T_{1,m})-I^{-1/2} J\big\| \big\|\widehat{\theta}(T_{\ell,k})-\theta^*_1\big\|\\
 &~~~~=o_P(1)+o_P(1)=o_P(1). 
\end{align*} 
\begin{flushright}
$\Box$
\end{flushright}
\subsubsection{Proof of Theorem \ref{th3}}
According to (\ref{Rejet_H0*}), it is enough to show that
\begin{equation*}
\sup_{m<k\leq [Tm]+1}\max_{\ell \in \Pi_{m,k}} \frac{\widehat{D}_{k,\ell}}{b((k-\ell)/m)}
\limiteloim
	 \sup_{1<t \leq T}\sup_{1<s<t} \frac{\| W_d(s)-sW_d(1)\|}{t b(s)}.
\end{equation*}
According to Lemma \ref{Lem_2}, to obtain the above convergence, it suffices to show that 
\begin{equation}\label{Eq1.th3}
\sup_{m<k\leq [Tm]+1}\max_{\ell \in \Pi_{m,k}} \frac{D_{k,\ell}}{b((k-\ell)/m)}
\limiteloim
	 \sup_{1<t \leq T}\sup_{1<s<t} \frac{\| W_d(s)-sW_d(1)\|}{t b(s)}.
\end{equation}
Let $k>m$ and $\ell \in \Pi_{m,k}$.  As $m \rightarrow +\infty$, we can proceed similarly as in (\ref{Eq_a_bis}) and (\ref{Eq_a_bisbis}) to show that  
\begin{align*}
&J\left(\widehat{\theta}(T_{1,m})-\theta^*_1\right)
 =
 \frac{1}{m}\frac{\partial}{\partial \theta} L(T_{1,m},\theta^*_1)+o_P\left(\frac{1}{\sqrt{m}}\right)  \text{ and}\\
 &J\left(\widehat{\theta}(T_{\ell,k})-\theta^*_1\right)
 =
 \frac{1}{k-\ell}\frac{\partial}{\partial \theta} L(T_{\ell,k},\theta^*_1)+o_P\left(\frac{1}{\sqrt{k-\ell}}\right).
\end{align*}
The subtraction of the two above equalities gives
  \[ 
  J (\widehat{\theta}(T_{\ell, k}) - \widehat{\theta}(T_{1, m}) ) = \frac{1}{k-\ell}  \Big( \frac{\partial}{\partial \theta} L(T_{\ell,k},\theta_1^*)
           - \frac{k-\ell}{m} \frac{\partial}{\partial \theta} L(T_{1,m},\theta_1^*)  \Big) +   o_P\left( \frac{1}{\sqrt{k-\ell}} + \frac{1}{\sqrt{m}}  \right)   .
    \]
This implies
\begin{equation*}
  \sqrt{m}\frac{k-\ell}{k}I^{-1/2}J (\widehat{\theta}(T_{\ell, k}) - \widehat{\theta}(T_{1, m}) ) = \frac{\sqrt{m}}{k}I^{-1/2}\Big( \frac{\partial}{\partial \theta} L(T_{\ell,k},\theta_1^*)-
           \frac{k-\ell}{m} \frac{\partial}{\partial \theta} L(T_{1,m},\theta_1^*)  \Big) +  o_P(1).
    \end{equation*}
Hence,
 \begin{align*}
  &\sup_{m<k\leq [Tm]+1}\max_{\ell \in \Pi_{m,k}} \frac{1}{ b((k-\ell)/m) }
  \Big\| \sqrt{m}\frac{k-\ell}{k}I^{-1/2}J  (\widehat{\theta}(T_{\ell,k}) - \widehat{\theta}(T_{1,m}) ) \\ 
 & \hspace{7cm} - \frac{ \sqrt{m}}{k} I^{-1/2} \Big( \frac{\partial}{\partial \theta} L(T_{\ell,k},\theta_1^*)    - \frac{k-\ell}{m} \frac{\partial}{\partial \theta} L(T_{1,m},\theta_1^*)  \Big)    \Big\| \\
  &\leq  
  \dfrac{1}{\inf_{s>0}b(s)}  \sup_{m<k\leq [Tm]+1}\max_{\ell \in \Pi_{m,k}}
  \Big\| \sqrt{m}\frac{k-\ell}{k}I^{-1/2}J (\widehat{\theta}(T_{\ell, k}) - \widehat{\theta}(T_{1, m}) ) \\
  & \hspace{6.5cm} - \frac{ \sqrt{m}}{k} I^{-1/2} \Big( \frac{\partial}{\partial \theta} L(T_{\ell,k},\theta_1^*)
           - \frac{k-\ell}{m} \frac{\partial}{\partial \theta} L(T_{1,m},\theta_1^*)  \Big)    \Big\| \\
 &~= 
 o_P(1).
 \end{align*}
Thus, to prove (\ref{Eq1.th3}),  we will show that
\begin{multline}\label{Eq2.th3}
\sup_{m<k\leq [Tm]+1}\max_{\ell \in \Pi_{m,k}} \frac{1}{ b((k-\ell)/m) }  \frac{ \sqrt{m}}{k}
  \Big\| I^{-1/2} \Big( \frac{\partial}{\partial \theta} L(T_{\ell,k},\theta_1^*)    - \frac{k-\ell}{m} \frac{\partial}{\partial \theta} L(T_{1,m},\theta_1^*)  \Big)    \Big\|\\
  \limiteloim
	 \sup_{1<t \leq T}\sup_{1<s<t} \frac{\| W_d(s)-sW_d(1)\|}{t b(s)}.
\end{multline}
~\\
Now, let us consider the following two cases.

 {\bf(1.)} Closed-end procedure.\\
 Let $1<T<\infty$. Define the set $S \coloneqq \left\{(t,s) \in [1,T]\times [1,T]/\, s<t\right\}$.
 According to Lemma \ref{Lem_1}, $\left(\frac{\partial}{\partial \theta} \ell_t(\theta_1^*),\mathcal{F}_{t}\right)_{t \in \mathbb{Z}}$ is a stationary ergodic martingale difference sequence with covariance matrix $I$.  
 Then, by Cramér-Wold device (see Billingsley (1968)), it holds that
\[
 \frac{1}{m} \sum_{i=[ms]+1}^{[mt]}\frac{\partial}{\partial \theta} \ell_i(\theta^*_1) \limiteloimds W_I(t-s),
\]
where $\limiteloimds$ means the weak convergence on the Skorohod space $\cal D(S)$ and $B_I$ is a $d$-dimensional Gaussian centered process such as $\E (B_I(s)B_I(\tau)') = \min(s,\tau )I$. Therefore,
\begin{align*}
& \frac{1}{\sqrt{m}} \Big(\sum_{i=[ms]+1}^{[mt]}\frac{\partial \ell_i(\theta^*_1)}{\partial \theta}-\frac{[mt]-[ms]}{m} \sum_{i=1}^{m}\frac{\partial \ell_i(\theta^*_1)}{\partial \theta}
   \Big)
   \limiteloimds B_I(t-s)-(t-s)B_I(1)\\
   &\text{and}~~~~~\\
   & \frac{1}{\sqrt{m}} I^{-1/2}\Big(\sum_{i=[ms]+1}^{[mt]}\frac{\partial \ell_i(\theta^*_1)}{\partial \theta}-\frac{[mt]-[ms]}{m} \sum_{i=1}^{m}\frac{\partial \ell_i(\theta^*_1)}{\partial \theta}
   \Big)
   \limiteloimds W_d(t-s)-(t-s)W_d(1).
   \end{align*}
Hence,
\begin{align*}
  & \sup_{m<k < Tm} \max_{ \ell \in \Pi_{m,k}} \frac{1}{ b((k-\ell)/m) } \frac{\sqrt{m}}{k}\Big \| I^{-1/2} \Big(\frac{\partial}{\partial \theta} L(T_{\ell,k},\theta_1^*)-
           \frac{k-\ell}{m} \frac{\partial}{\partial \theta} L(T_{1,m},\theta_1^*)\Big)  \Big \|\\
   & ~~=
   \sup_{m<k < Tm} \max_{ \ell \in \Pi_{m,k}} \frac{1}{ b((k-\ell)/m) } \frac{\sqrt{m}}{k}
   \Big \| I^{-1/2} \Big(
   \sum_{i=\ell}^{k}\frac{\partial \ell_i(\theta^*_1)}{\partial \theta}-\frac{k-\ell}{m} \sum_{i=1}^{m}\frac{\partial \ell_i(\theta^*_1)}{\partial \theta}
   \Big)
   \Big \|\\
   & ~~=
   \sup_{t \in \{1,1+\frac{1}{m},\cdots,T \}} \max_{ s \in \{1,1-\frac{v'_m}{m}, 2-\frac{v'_m}{m},\cdots,t-\frac{v'_m}{m}\}}
   \Bigg[
   \frac{1}{b(([mt]-[ms])/m)} \frac{m}{[mt]}\\
   & \hspace{7cm}
   \times\Big \|\frac{1}{\sqrt{m}} I^{-1/2}\Big(\sum_{i=[ms]+1}^{[mt]}\frac{\partial \ell_i(\theta^*_1)}{\partial \theta}-\frac{[mt]-[ms]}{m} \sum_{i=1}^{m}\frac{\partial \ell_i(\theta^*_1)}{\partial \theta}
   \Big)\Big \|
   \Bigg]\\
   & \hspace{0.3cm} 
   \limiteloim
           \sup_{1<t < T}\sup_{1<s<t} \frac{\| W_d(t-s)-(t-s)W_d(1)\|}{t b(t-s)}
          \stackrel{\cal D}{=}
           \sup_{1<t < T}\sup_{1<s<t} \frac{\| W_d(s)-sW_d(1)\|}{t b(s)}.
\end{align*}
Thus, (\ref{Eq2.th3}) follows; which ends the proof in the case of the closed-end procedure.\\

{\bf (2.)} Open-end procedure.\\
According to (\ref{Eq2.th3}) and  {\bf (1.)}, it suffices to show that the limit distribution (as $m,T \rightarrow \infty $) of
\[
\sup_{k > [Tm]}  \max_{ \ell \in \Pi_{m,k}} \frac{1}{ b((k-\ell)/m) } \frac{\sqrt{m}}{k}\Big \| I^{-1/2}\Big(\frac{\partial}{\partial \theta} L(T_{\ell,k},\theta_1^*)-
           \frac{k-\ell}{m} \frac{\partial}{\partial \theta} L(T_{1,m},\theta_1^*)\Big)  \Big \|
\]
exists and is equal to the limit distribution (as $T \rightarrow \infty $) of
\[
\sup_{t > T}\sup_{1<s<t} \frac{\| W_d(s)-sW_d(1)\|}{tb(s)}.
\]
Let $k > [Tm]$. For some $\ell_k \in \Pi_{m,k}$, we have
\[
\max_{ \ell \in \Pi_{m,k}} \frac{1}{ b((k-\ell)/m) } \frac{\sqrt{m}}{k}\Big \| \frac{\partial}{\partial \theta} L(T_{\ell,k},\theta_1^*) \Big \| 
= \frac{1}{ b((k-\ell_k)/m) } \frac{\sqrt{m}}{k}  \Big \|\sum_{i=\ell_k}^{k}\frac{\partial}{\partial \theta}\ell_i(\theta^*_1)\Big \|.
\]
From the Hájek-Rényi-Chow inequality (see Chow (1960)), we can get
\begin{equation}\label{Eq3.th3}
\forall \varepsilon >0,~~~
\lim_{T \rightarrow \infty} \limsup_{m \rightarrow \infty} \prob \Big( \sup_{k>Tm} \frac{1}{ b((k-\ell_k)/m) }\frac{\sqrt{m}}{k}  \Big \|\sum_{i=\ell_k}^{k}\frac{\partial}{\partial \theta}\ell_i(\theta^*_1)\Big \| > \varepsilon \Big)=0.
\end{equation}
Moreover, since the function $b(\cdot)$ is non-increasing, for any $m, T > 1$, we have 
\begin{align}\label{Eq4.th3}
\sup_{k > Tm}  \max_{ \ell \in \Pi_{m,k}} \frac{1}{ b((k-\ell)/m) } \frac{\sqrt{m}}{k}\Big \| \frac{k-\ell}{m}\frac{\partial}{\partial \theta} L(T_{1,m},\theta_1^*) \Big \| 
&= 
 \Big \| \frac{1}{\sqrt{m}}   \sum_{i=1}^{m}\frac{\partial}{\partial \theta}\ell_i(\theta^*_1)\Big \|
 \cdot 
 \sup_{k > Tm}  \max_{ \ell \in \Pi_{m,k}} \frac{1}{ b((k-\ell)/m) } \frac{k-\ell}{k} \nonumber\\
 &=
 \Big \| \frac{1}{\sqrt{m}}   \sum_{i=1}^{m}\frac{\partial}{\partial \theta}\ell_i(\theta^*_1)\Big \|
 \cdot 
 \sup_{k > Tm} \frac{1}{ b((k-v'_m)/m) } \frac{k-v'_m}{k}\nonumber\\
 &=
  \frac{1}{\inf_{s>0}b(s)}   \Big \| \frac{1}{\sqrt{m}}   \sum_{i=1}^{m}\frac{\partial}{\partial \theta}\ell_i(\theta^*_1)\Big \|\nonumber\\
 &=
  \limiteloim \frac{1}{\inf_{s>0}b(s)}  \|B_I(1)  \|,
\end{align}
by using again the Cramèr-Wold device and the central limit theorem applied to the martingale difference sequence $\left(\frac{\partial}{\partial \theta} \ell_t(\theta_1^*),\mathcal{F}_{t}\right)_{t \in \mathbb{Z}}$.
It comes from (\ref{Eq3.th3}) and (\ref{Eq4.th3}) that
\begin{multline}\label{Eq5.th3}
\sup_{k > [Tm]}  \max_{ \ell \in \Pi_{m,k}} \frac{1}{ b((k-\ell)/m) } \frac{\sqrt{m}}{k}\Big \| I^{-1/2}\Big(\frac{\partial}{\partial \theta} L(T_{\ell,k},\theta_1^*)-
           \frac{k-\ell}{m} \frac{\partial}{\partial \theta} L(T_{1,m},\theta_1^*)\Big)  \Big \|\\
          \overset{\mathcal{D}}{\underset{m, T \rightarrow \infty}{\longrightarrow} } \frac{1}{\inf_{s>0}b(s)}  \|W_d(1)  \|.
\end{multline}
Furthermore, from the proof of Lemma 6.3 of Bardet and Kengne (2014), we get
\[
\sup_{t> T}\sup_{1<s<t} \frac{\| B_I(s)-sB_I(1)\|}{t b(s)} \overset{\mathcal{D}}{\underset{ T \rightarrow \infty}{\longrightarrow} }  \frac{1}{\inf_{s>0} b(s)} \| B_I(1)\| ~~ \text{as} ~T \rightarrow \infty.
\]
Therefore,
\begin{equation}\label{Eq6.th3}
\sup_{t> T}\sup_{1<s<t} \frac{\| W_d(s)-sW_d(1)\|}{t b(s)} \overset{\mathcal{D}}{\underset{T \rightarrow \infty}{\longrightarrow} }  \frac{1}{\inf_{s>0} b(s)} \| W_d(1)\| ~~ \text{as} ~T \rightarrow \infty.
\end{equation}
The relations (\ref{Eq5.th3}) and (\ref{Eq6.th3}) complete the proof in the case of the open-end procedure.
\begin{flushright}
$\blacksquare$
\end{flushright}

\subsubsection{Proof of Theorem \ref{th4}}
Denote $k_m= k^* + m^\delta$ for $\delta \in (1/2,1)$.
 For $m$ large enough, we have $v'_m <m^\delta $ and thus $k^* \leq k_m-v'_m$.
 Moreover since $k^*(m)=[T^*m]$ for some $T^*<T$, for $m$ large enough and for both the open-end and closed-end
procedures, $k_m= k^* + m^\delta< T^*m + m^\delta<Tm$. 
Hence, $k_m$ is between $m$ and $[Tm]+1$ and $k^* \in \Pi_{m,k_m}$ for $m$ large enough.
 Therefore, according to Assumption $\bf B_*$, we can find a constant $C>0$ such that
 \begin{align} \label{Eq7._th4} 
  \underset{ \ell \in \Pi_{m,k_m}}{\max} \frac{\widehat{D}_{k_m,\ell}}{ b((k_m-\ell)/m)}
    &=  
            \max_{\ell \in \Pi_{m,k_m}} \frac{1}{ b((k_m-\ell)/m) }
            \sqrt{m} \frac{k_m-\ell}{k_m} \big\| \widehat I(T_{1,m})^{-1/2} \widehat J(T_{1,m})  \big(\widehat{\theta}(T_{\ell,k_m}) - \widehat{\theta}(T_{1,m})\big)\big\|  \nonumber\\
 \nonumber 
 & \geq  
 \dfrac{1}{ b((k_m-k^*)/m) }
            \sqrt{m} \frac{k_m-k^*}{k_m} \big\| \widehat I(T_{1,m})^{-1/2} \widehat J(T_{1,m}) \big(\widehat{\theta}(T_{k^*,k_m}) - \widehat{\theta}(T_{1,m})\big)\big\|  \nonumber \\
 \nonumber 
 & \geq 
  C \sqrt{m}  \dfrac{m^\delta}{[T^* m] + m^\delta} \big\| \widehat I(T_{1,m})^{-1/2} \widehat J(T_{1,m}) \big(\widehat{\theta}(T_{k^*,k_m}) - \widehat{\theta}(T_{1,m})\big)\big\| \nonumber\\
& \geq  C m^{\delta - 1/2}  \big\| \widehat I(T_{1,m})^{-1/2} \widehat J(T_{1,m})  (\widehat{\theta}(T_{k^*,k_m}) - \widehat{\theta}(T_{1,m}))\big\| .
   \end{align}
   Moreover, from \cite{Francq2016} , we get 
    \[ 
    \widehat I(T_{1,m})  \limitepsm I , ~~~  \widehat J(T_{1,m}) \limitepsm  J, ~~~ 
     \widehat{\theta}(T_{1,m}) \limitepsm \theta^*_1  ~~~\text{and} ~~~
     \widehat{\theta}(T_{k^*,k_m}) \limitepsm \theta^*_2 .
   \]
  Thus, since $I$ and $J$ are symmetric positive definite,  $\theta^*_1 \neq \theta^*_2$ and $  \delta > 1/2 $, it comes from (\ref{Eq7._th4}) that 
 \[ \underset{ \ell \in \Pi_{m,k_m}}{\max} \frac{\widehat{D}_{k_m,\ell}}{ b((k_m-\ell)/m)}  \limitepsm  \infty  . \]
 \begin{flushright}
$\blacksquare$
\end{flushright}



 \end{document}